\documentclass[11pt,a4paper,reqno]{amsart}
\usepackage{graphicx}
\usepackage[linktocpage=true,colorlinks,citecolor=blue,linkcolor=blue,urlcolor=blue]{hyperref}
\usepackage{amssymb}
\usepackage{cite}
\usepackage{amsmath}
\usepackage{latexsym}
\usepackage{amscd}
\usepackage{amsthm}
\usepackage{mathrsfs}
\usepackage{url}

\newtheorem{thm}{Theorem}[section]
\newtheorem{corr}[thm]{Corollary}
\newtheorem{lem}[thm]{Lemma}
\newtheorem{prop}[thm]{Proposition}

\theoremstyle{definition}

\theoremstyle{remark}
\newtheorem{rem}[thm]{Remark}

\numberwithin{equation}{section}
\def\R{\mathbb R}

\def\C{\mathbb C}
\def\P{\mathbb P}

\def\ra{\rightarrow}
\def\pt{\partial}
\DeclareMathOperator\GG{G}
\DeclareMathOperator\GL{GL}
\DeclareMathOperator\SO{SO}
\DeclareMathOperator\Diff{Diff}

\begin{document}
\title[Stability of torsion-free \texorpdfstring{$\GG_2$}{G2} structures]{Stability of torsion-free \texorpdfstring{{\boldmath $\GG_2$}}{G2} structures along the Laplacian flow}
\author{Jason D. Lotay}
\author{Yong Wei}
\address{Department of Mathematics, University College London, Gower Street, London, WC1E 6BT, United Kingdom}
\email{j.lotay@ucl.ac.uk, yong.wei@ucl.ac.uk}
\thanks{This research was supported by EPSRC grant EP/K010980/1.}

\begin{abstract}
We prove that torsion-free $\GG_2$ structures are (weakly) dynamically stable along the Laplacian flow for closed $\GG_2$ structures. More precisely, given a torsion-free $\GG_2$ structure $\bar{\varphi}$ on a compact $7$-manifold $M$, the Laplacian flow with initial value
 in $[\bar{\varphi}]$, sufficiently close to $\bar{\varphi}$, will converge to a point in the $\Diff^0\!(M)$-orbit of $\bar{\varphi}$.  We deduce, from fundamental work of Joyce \cite{joyce96-1}, that the Laplacian flow starting 
  at any closed $\GG_2$ structure with sufficiently small torsion will exist for all time and converge to a torsion-free $\GG_2$ structure.
\end{abstract}

\maketitle

\section{Introduction}

In this article we focus on the Laplacian flow for closed $\GG_2$ structures, whose stationary points are torsion-free $\GG_2$ structures on a compact 7-manifold.  This geometric flow thus gives a potential means for studying the challenging problem of existence of 
 Ricci-flat metrics with exceptional holonomy $\GG_2$.  

A desirable property of a 
geometric flow is (weak) dynamical stability of its stationary points: if one starts the flow close to a stationary point then it will 
converge to a (possibly different) stationary point.  Such a property gives one hope to obtain further long-time existence and convergence results.   
We prove that weak dynamical stability holds for torsion-free $\GG_2$ structures along the Laplacian flow, building on our work 
in \cite{Lotay-Wei}, thus confirming an expectation of Bryant \cite{bryant2005}  

As a consequence, we deduce from Joyce's fundamental perturbative existence theory for torsion-free $\GG_2$ structures \cite{joyce96-1}
that the Laplacian flow starting at a closed $\GG_2$ structure with torsion sufficiently close to $0$, in a suitable sense, will exist for all time and 
converge to a torsion-free $\GG_2$ structure.  This  provides key evidence that the Laplacian flow may indeed provide a useful tool for finding
holonomy $\GG_2$ metrics.

\medskip

A $\GG_2$ structure on a 7-manifold $M$ is
defined by a positive $3$-form $\varphi\in\Omega^3_+(M)$, where
the positivity is a natural nondegeneracy condition.  A positive 3-form exists if and only if $M$ is orientable and spinnable.  To each $\varphi\in\Omega^3_+(M)$, one can associate a unique Riemannian metric $g_{\varphi}$ and an orientation on $M$. 
 If $\nabla_{\varphi}$ is the Levi-Civita connection of $g_{\varphi}$, we can interpret $\nabla_{\varphi}\varphi$ as the torsion of the $\GG_2$ structure $\varphi$, so if $\nabla_{\varphi}\varphi=0$ we say $\varphi$ is torsion-free and that $(M,\varphi)$ is a $\GG_2$ manifold.  Equivalently, if $*_{\varphi}$ is the Hodge star operator determined by $g_{\varphi}$ and the orientation, then $\varphi$ is torsion-free if and only if $d\varphi=d\! *_\varphi\!\varphi=0$.

 The key property of torsion-free $\GG_2$ structures 
 is that the holonomy group $\textrm{Hol}(g_{\varphi})\subseteq \GG_2$, and thus $\GG_2$ manifolds $(M,\varphi)$ are Ricci-flat.  Moreover, one can characterise the compact $\GG_2$ manifolds with  $\textrm{Hol}(g_{\varphi})= \GG_2$ as those with finite fundamental group. Thus understanding torsion-free $\GG_2$ structures is crucial for constructing Riemannian manifolds with holonomy  $\GG_2$.

Bryant \cite{bryant1987} used the theory of exterior differential systems to demonstrate the local existence of many metrics with holonomy $\GG_2$, and then Bryant-Salamon \cite{bryant-salamon} constructed complete non-compact manifolds with holonomy $\GG_2$, which are the spinor bundle of $\mathcal{S}^3$ and the bundles of anti-self-dual 2-forms on $\mathcal{S}^4$ and $\C\P^2$. In \cite{joyce96-1}, Joyce constructed the first examples of compact $7$-manifolds with holonomy $\GG_2$ and many further
 compact examples have now been constructed \cite{Kov, CHNP}.

\medskip

 Since Hamilton \cite{ha82} introduced the Ricci flow in 1982, geometric flows have been an important tool in studying geometric structures on manifolds. In 1992, Bryant (see \cite{bryant2005}) proposed the Laplacian flow for closed $\GG_2$ structures 
\begin{equation}\label{Lap-flow-def}
  \left\{\begin{array}{rcl}
         \frac{\pt}{\pt t}\varphi &
         = &\Delta_{\varphi}\varphi,\\[2pt]
           d\varphi &
           = & 0, \\[2pt]
           \varphi(0) &
           = & \varphi_0,
         \end{array}\right.
\end{equation}
where $\Delta_{\varphi}=dd^*_{\varphi}+d^*_{\varphi}d$ is the Hodge Laplacian with respect to $g_{\varphi}$ and $\varphi_0$ is an initial closed $\GG_2$ structure.  The stationary
points of the flow are harmonic $\varphi$, which on a compact
manifold are precisely the torsion-free $\GG_2$ structures, so the Laplacian flow provides a tool for studying the existence of torsion-free $\GG_2$ structures 
on a manifold admitting closed $\GG_2$ structures.   We remark that there are other proposed flows which also
have torsion-free $\GG_2$ structures as stationary points (e.g.~\cite{Grig, Kar-solit, weiss-witt}).

The goal is to understand the long-time behavior of the Laplacian flow on compact manifolds $M$; specifically, to understand conditions under which the flow will converges to a torsion-free $\GG_2$ structure. A fundamental result of Joyce \cite{joyce96-1}  states that a closed $\GG_2$ structure with sufficiently small torsion (in a suitable sense) can be perturbed to a torsion-free $\GG_2$ structure in its cohomology class.
A reasonable first conjecture is therefore that a dynamic version of
Joyce's result holds: namely, that if the initial $\GG_2$ structure $\varphi_0$ on $M$ is closed and has small torsion, then the Laplacian flow \eqref{Lap-flow-def} exists for all time and converges to a torsion-free $\GG_2$ structure.  We show that such a conjecture is true in Corollary \ref{main-cor} below.

As an essential first step in studying \eqref{Lap-flow-def}, a short-time existence result was claimed in \cite{bryant2005} and  proved in detail in \cite{bryant-xu2011}.
\begin{thm}
\label{thm-bryant-xu}
 For a compact $7$-manifold $M$ with closed $\GG_2$
 structure $\varphi_0$, the initial value problem \eqref{Lap-flow-def} has a unique solution for a short time $t\in[0,\epsilon)$.
\end{thm}
\begin{rem}
The existence time $\epsilon$ depends on the initial value $\varphi_0$.
\end{rem}
To prove Theorem \ref{thm-bryant-xu}, Bryant--Xu \cite{bryant-xu2011} modified \eqref{Lap-flow-def} by a term $\mathcal{L}_{V(\varphi)}\varphi$ for some vector field $V(\varphi)$ and considered a
``Laplacian--DeTurck flow'':
\begin{equation}\label{deturck-flow}
\left\{\begin{array}{ccl}
        \frac{\pt}{\pt t}\varphi&=&\Delta_{\varphi}\varphi+\mathcal{L}_{V(\varphi)}\varphi,\\[2pt]
         d\varphi & = & 0, \\[2pt]
           \varphi(0) & = &\varphi_0,
         \end{array}\right.
\end{equation}
where $\mathcal{L}_{V(\varphi)}\varphi$ is the Lie derivative of $\varphi$ in the direction $V(\varphi)$. Linearizing  \eqref{deturck-flow},  one sees that it is parabolic in the direction of closed forms: this is not a typical type of parabolicity, and so standard parabolic theory does not obviously apply.  The conclusion therefore follows by applying  DeTurck's trick \cite{deturck1983} and the Nash--Moser inverse function theorem (c.f.~\cite{ha82,ha82-nash-moser}).
This result is particularly surprising since \eqref{Lap-flow-def} uses the Hodge Laplacian and so appears to have the wrong sign to be parabolic, but it turns out that by gauge-fixing as in \eqref{deturck-flow} one does indeed obtain a parabolic flow (at least in the required directions given by closed forms).

Hitchin \cite{hitchin2000} (see also \cite{bryant-xu2011}) showed that the Laplacian flow is related to a volume functional, which we now briefly describe.
Suppose that $\bar{\varphi}\in\Omega^3_+(M)$ is a closed $\GG_2$ structure on $M$. Let
\begin{equation*}
  [\bar{\varphi}]_+=\{\bar{\varphi}+d\beta\in\Omega^3_+(M)\,:\,\beta\in\Omega^2(M)\}
\end{equation*}
be the open subset of the cohomology class $[\bar{\varphi}]$ that consists of positive $3$-forms. The volume functional $\mathcal{H}:[\bar{\varphi}]_+\ra\R^+$ is defined by
\begin{equation*}
  \mathcal{H}(\varphi)=\frac 17\int_M\varphi\wedge *_{\varphi}\varphi=\int_M*_{\varphi}1. 
\end{equation*}
Then $\varphi\in[\bar{\varphi}]_+$ is a critical point of $\mathcal{H}$ if and only if $d *_{\varphi}\!\varphi=0$, 
i.e.~$\varphi$ is torsion-free.  Moreover, the Laplacian flow can be viewed as the gradient flow for $\mathcal{H}$, 
with respect to a non-standard $L^2$-type metric on $[\bar{\varphi}]_+$ (see e.g.~\cite{bryant-xu2011}).

 If $\bar{\varphi}$ is torsion-free then by calculating the second variation of the functional $\mathcal{H}$ at $\bar{\varphi}$, one sees that the orbit $\Diff^0(M)\cdot\bar{\varphi}$ of $\bar{\varphi}$ under the diffeomorphisms of $M$ that are isotopic to identity is a local maximum of $\mathcal{H}$ on the moduli space $[\bar{\varphi}]_+/\Diff^0(M)$ (see \cite{hitchin2000,bryant2005}). Consequently,  Bryant \cite{bryant2005} expected that for $\varphi_0\in [\bar{\varphi}]_+$ sufficiently close to $\bar{\varphi}$, 
the Laplacian flow \eqref{Lap-flow-def} starting at $\varphi_0$ will converge to a point in the orbit $\Diff^0(M)\cdot\bar{\varphi}$.
Our main result confirms this expectation.

\begin{thm}\label{thm-main}
Let $\bar{\varphi}$ be a torsion-free $\GG_2$ structure  
on a compact $7$-manifold $M$. There is a neighborhood
$\mathcal{U}$ of $\bar{\varphi}$ such that for any $\varphi_0\in [\bar{\varphi}]_+\cap \mathcal{U}$, the Laplacian flow \eqref{Lap-flow-def} 
with initial value $\varphi_0$ exists for all $t\in [0,\infty)$ and converges smoothly to $\varphi_{\infty}\in \Diff^{0}(M)\cdot\bar{\varphi}$ as $t\ra\infty$.

In other words,  torsion-free $\GG_2$ structures are (weakly) dynamically stable along the Laplacian flow for closed $\GG_2$ structures.
\end{thm}

The proof of Theorem \ref{thm-main} is inspired by the proof of an analogous result in Ricci flow: Ricci-flat metrics are dynamically stable along the Ricci
flow.  The idea is to combine arguments for the Ricci flow case \cite{has, sesum2004, sesum2006} with the particulars of the geometry of closed $\GG_2$
structures and new higher order
estimates for the Laplacian flow derived by the authors \cite{Lotay-Wei}.  For a precise definition of the neighbourhood $\mathcal{U}$, which is defined by
norms involving derivatives at least up to order $9$, see \S\ref{sec-proof of thm}.  

If we assume Joyce's existence result for torsion-free $\GG_2$ structures \cite{joyce96-1}, then Theorem \ref{thm-main} has the following immediate corollary.

\begin{corr}\label{main-cor}  Let $\varphi_0$ be a closed $\GG_2$ structure on a compact 7-manifold $M$.
There exists a $C^{9}$-open neighbourhood $\mathcal{U}$ of $0$ in $\Omega^3(M)$ such that if $d^*_{\varphi_0}\varphi_0=d^*_{\varphi_0}\gamma$ for some
 $\gamma\in\mathcal{U}$, then the Laplacian flow \eqref{Lap-flow-def} with initial value $\varphi_0$ exists for all time and converges to
 a torsion-free $\GG_2$ structure.
\end{corr}

\noindent The theory in \cite{joyce96-1} states that if we control the $C^0$ and $L^2$-norms of $\gamma$ and the $L^{14}$-norm of
$d^*_{\varphi_0}\gamma=d^*_{\varphi_0}\varphi_0$,  we can deform $\varphi_0$ in its cohomology class to a unique $C^0$-close torsion-free $\GG_2$ structure
$\bar{\varphi}$.  By choosing a neighbourhood $\mathcal{U}$ appropriately, controlling derivatives up to at least order $9$, we can ensure that we can apply both
the theory of \cite{joyce96-1}
and Theorem \ref{thm-main}, and thus deduce the corollary.  The neighbourhood $\mathcal{U}$ given by Corollary \ref{main-cor} is not optimal, and
one would like to able to prove this result directly  using the Laplacian flow with optimal conditions and without recourse to \cite{joyce96-1}, but 
nevertheless, Corollary \ref{main-cor} gives significant evidence that the Laplacian flow will play an important role in understanding the problem of existence of torsion-free $\GG_2$ structures on 7-manifolds admitting closed $\GG_2$ structures. 

Our results also motivate us to study an approach to the following  problem\footnote{The authors are grateful to Thomas Walpuski for pointing this out.}.  The work of Joyce \cite{joyce96-1} 
shows that the natural map  from the moduli space $\mathcal{M}$ of torsion-free $\GG_2$ structures (that is, torsion-free
$\GG_2$ structures up to the action of $\Diff^0(M)$) to $H^3(M)$ given by 
$\Diff^0(M)\cdot\bar{\varphi}\mapsto [\bar{\varphi}]$ is locally 
injective, but the question of whether this map is globally injective, raised by Joyce (c.f.~\cite{joyce2000}), is still open.  Equivalently, it can be phrased as asking whether the Lagrangian immersion 
of $\mathcal{M}$ into $T^*H^3(M)\cong H^3(M)\times H^4(M)$ given by $$\Diff^0(M)\cdot\bar{\varphi}\mapsto ([\bar{\varphi}],[*_{\bar{\varphi}}\bar{\varphi}])$$
is an embedding.
  Suppose we have two torsion-free $\GG_2$ structures $\bar{\varphi}_0$ and $\bar{\varphi}_1$ which lie in the same cohomology class, 
so we can write $\bar{\varphi}_1=\bar{\varphi}_0+ d\eta$ for some 2-form $\eta$. We would like to see whether $\bar{\varphi}_1\in\Diff^0(M)\cdot\bar{\varphi}_0$.   
  By our main theorem (Theorem \ref{thm-main}) we know that the Laplacian flow starting at $\varphi_0(s)=\bar{\varphi}_0+sd\eta$ (which is closed) will exist for all time and converge to $\phi_s^*\bar{\varphi}_0$ for some $\phi_s\in\Diff^0(M)$ when $s$ is sufficiently small.  Similarly, the Laplacian flow starting 
  at $\varphi_0(s)$ for $s$ near $1$ will also exist for all time and now converge to $\Phi_s^*\bar{\varphi}_1$ for some $\Phi_s\in\Diff^0(M)$.  
  The aim would be to study long-time existence and convergence of the flow starting at any $\varphi_0(s)$ and, if this occurs, use 
Joyce's local injectivity result to show the $\GG_2$ moduli space $\mathcal{M}$ embeds into $T^*H^3(M)$, or otherwise find obstructions.

Generally, one cannot expect the Laplacian flow \eqref{Lap-flow-def} will converge to a torsion-free $\GG_2$ structure, even if it has long-time existence. There are compact 7-manifolds with closed
$\GG_2$ structures that cannot admit holonomy $\GG_2$ metrics  for topological reasons (c.f.~\cite{fernandez1987-1,fernandez1987-2}), and Bryant \cite{bryant2005} showed that the Laplacian flow  starting with a particular one of these examples will exist for all time but it does not converge; for instance,  the volume of the associated metrics will increase without bound.
However, we note that recently Fern\'{a}ndez--Fino--Manero \cite{Fern-Fino-M} constructed some non-compact solutions for the Laplacian flow which exist for all time and do converge, but to a flat $\GG_2$ structure, and Bryant \cite{Bryantsol} has informed the first author of non-trivial examples of non-compact steady solitons for the Laplacian flow (i.e.~solutions which exist for all time and only vary within their diffeomorphism orbit).

\medskip

We now briefly summarise the contents of this article.  
In $\S$\ref{sec-prelim} we recall some basic properties of closed $\GG_2$ structures and the Laplacian flow, mainly from \cite{bryant2005} but also with some new observations from \cite{Lotay-Wei}.  For the remainder of the article we assume that $(M,\bar{\varphi})$ is a compact $\GG_2$ manifold.
In $\S$\ref{sec-estimate}, we study the Laplacian--DeTurck flow \eqref{deturck-flow} for $\varphi_0\in[\bar{\varphi}]_+$ and show that
if $\varphi_0$ is sufficiently close to $\bar{\varphi}$ then the flow will exist at least up to time $1$ and remain close to $\bar{\varphi}$.
In $\S$\ref{sec:t>1} we show that the Laplacian--DeTurck flow exponentially decays in $L^2$ to $\bar{\varphi}$ and that we
can use Sobolev estimates to show that as long as the flow exists, it
is uniformly controlled by the initial error $\varphi_0-\bar{\varphi}$. In \S \ref{sec-proof of thm}, we first prove dynamical stability of torsion-free $\GG_2$ 
structures under the Laplacian--DeTurck flow \eqref{deturck-flow} for closed $\GG_2$ structures, and then use this to prove our main result, Theorem \ref{thm-main}.

\section{Preliminaries}\label{sec-prelim}

In this section, we collect some facts and results on closed $\GG_2$ structures and the Laplacian flow. For more detail, we refer the reader to \cite{bryant2005, bryant-xu2011, joyce2000, Lotay-Wei, salamon1989}.

Let $\{e_1,e_2,\cdots,e_7\}$ denote the standard basis of $\R^7$ and
let  $\{e^1,e^2,\cdots,e^7\}$ be its dual basis. Write $e^{ijk}=e^i\wedge e^j\wedge e^k$ for simplicity and define a $3$-form $\phi$ by
\begin{equation*}
  \phi=e^{123}+e^{145}+e^{167}+e^{246}-e^{257}-e^{347}-e^{356}.
\end{equation*}
The subgroup of $\GL(7,\R)$ fixing $\phi$ is the exceptional Lie group $\GG_2$, which is a compact, connected, simple Lie subgroup of $\SO(7)$ of dimension $14$.
Note that $\GG_2$ acts irreducibly on $\R^7$ and preserves the metric and orientation for which $\{e_1,e_2,\cdots,e_7\}$ is an oriented orthonormal basis. If $*_{\phi}$ denotes the Hodge star on $\R^7$, then $\GG_2$ also preserves the $4$-form
\begin{equation*}
  *_{\phi}\phi=e^{4567}+e^{2367}+e^{2345}+e^{1357}-e^{1346}-e^{1256}-e^{1247}.
\end{equation*}

Let $M$ be a $7$-manifold. For $x\in M$ we let
\begin{equation*}
  \Lambda^3_+(M)_x=\{\varphi_x\in\Lambda^3T_x^*M|\exists u\in\textrm{Hom}(T_xM,\R^7), u^*\phi=\varphi_x\},
\end{equation*}
and thus obtain a bundle $\Lambda^3_+(M)=\sqcup_x \Lambda^3_+(M)_x$, which is an open subbundle of $\Lambda^3T^*M$, with fibre $\GL(7,\mathbb{R})/{\GG_2}$. We call a section $\varphi$ of $\Lambda^3_+(M)$ a positive  $3$-form on $M$ and we denote the space of sections of $\Lambda^3_+(M)$ by $\Omega^3_+(M)$. There is a 1-1 correspondence between $\GG_2$ structures and positive  $3$-forms, because given $\varphi\in\Omega^3_+(M)$, the subbundle of the frame bundle 
whose fibre at $x$ consists of $u\in\textrm{Hom}(T_xM,\R^7)$
such that $u^*\phi=\varphi_x$  defines a principal subbundle with fibre $\GG_2$, i.e.~a $\GG_2$ structure on $M$. Thus we usually call $\varphi\in\Omega^3_+(M)$ a $\GG_2$ structure on $M$.  Note that the existence of $\GG_2$ structures is a purely topological fact, namely that $M$ is oriented and spin.

The $\GG_2$ structure induces a splitting of the bundle of $k$-forms   ($2\leq k\leq 5$)
into direct summands, which we denote by $\Lambda^k_l(T^*M,\varphi)$ so that $l$ indicates the rank of the bundle.   If we let $\Omega^k_l(M)$ be the space of sections of $\Lambda^k_l(T^*M,\varphi)$ then we have that
\begin{align*}
  \Omega^2(M)= & \Omega^2_7(M)\oplus\Omega^2_{14}(M), \\
   \Omega^3(M)=& \Omega^3_1(M)\oplus\Omega^3_7(M)\oplus \Omega^3_{27}(M),
\end{align*}
where, if we let $\psi=*_{\varphi}\varphi$,
\begin{align*}
  \Omega^2_7(M) =& \{\beta\in\Omega^2(M)|\beta\wedge\varphi=2*_{\varphi}\beta\} =\{X\lrcorner\varphi|X\in\Gamma(TM)\},\\
  \Omega^2_{14}(M) =& \{\beta\in\Omega^2(M)|\beta\wedge\varphi=-*_{\varphi}\beta\}=\{\beta\in\Omega^2(M)|\beta\wedge\psi=0\}
\end{align*}
and
\begin{align*}
  \Omega^3_1(M)=&\{f\varphi|f\in C^{\infty}(M)\},\\
  \Omega^3_7(M)=&\{X\lrcorner\psi|X\in \Gamma(TM)\},\\
  \Omega^3_{27}(M)=&\{\gamma\in\Omega^3(M)|\gamma\wedge\varphi=0=\gamma\wedge\psi\}.
\end{align*}
By Hodge duality, we have a similar decomposition of $\Omega^4(M)$ and $\Omega^5(M)$.

To study the Laplacian flow, it is convenient to write the relevant quantities for study in local coordinates.
We write a $k$-form $\alpha$ locally using summation convention and  index notation as follows:
\begin{equation*}
  \alpha=\frac 1{k!}\alpha_{i_1i_2\cdots i_k}dx^{i_1}\wedge\cdots\wedge dx^{i_k},
\end{equation*}
where $\{x^1,\cdots, x^7\}$ give  local coordinates  on $M$.

We define the operator $i_{\varphi}: S^2T^*M\ra \Lambda^3T^*M$ as in \cite{bryant2005} (up to a factor of $\frac{1}{2}$) by
\begin{align*}
  i_{\varphi}(h)=&\frac 12h^l_i\varphi_{ljk}dx^i\wedge dx^j\wedge dx^k,
\end{align*}
for $h=h_{ij}dx^idx^j$ locally. Then $\Lambda^3_{27}(T^*M,\varphi)=i_{\varphi}(S^2_0T^*M)$, where $S^2_0T^*M$ denotes the trace-free symmetric $2$-tensors on $M$, and $i_{\varphi}(g)=3\varphi$. We also have the  map $j_{\varphi}$ on 3-forms:
\begin{equation}\label{j-varphi-def}
  j_{\varphi}(\gamma)(u,v)=*_{\varphi}((u\lrcorner\varphi)\wedge (v\lrcorner\varphi)\wedge \gamma)
\end{equation}
which is an isomorphism between $\Omega^3_1(M)\oplus \Omega^3_{27}(M)$ and the symmetric $2$-tensors on $M$  that provides an inverse of $i_{\varphi}$. Then  we have
\begin{equation*}
  j_{\varphi}(i_{\varphi}(h))=4h+2tr_g(h)g.
\end{equation*}
for any $h\in S^2(M)$, and $j_{\varphi}(\varphi)=6g$.

Given any $\GG_2$ structure $\varphi\in\Omega^3_+(M)$, there exist unique differential forms $\tau_0\in\Omega^0(M), \tau_1\in\Omega^1(M), \tau_{2}\in\Omega^2_{14}(M)$ and $\tau_{3}\in\Omega^3_{27}(M)$ such that $d\varphi$ and $d\psi$ can be expressed as follows:
\begin{align}
  d\varphi =& \tau_0\psi+3\tau_1\wedge\varphi+*_{\varphi}\tau_{3},\label{dvarphi} 
\\
  d\psi= & 4\tau_1\wedge\psi+\tau_{2}\wedge\varphi.\label{dpsi}
\end{align}
We call the four forms $\{\tau_0,\tau_1,\tau_2,\tau_3\}$ the intrinsic  torsion forms of a $\GG_2$ structure.  The full torsion tensor is a $2$-tensor $T_{ij}$ satisfying
\begin{equation}\label{nabla-var}
  \nabla_i\varphi_{jkl}=T_i^{\,\,m}\psi_{mjkl},
\end{equation}
\begin{equation}\label{T-def}
  T_i^{\,\,j}=\frac 1{24}\nabla_i\varphi_{lmn}\psi^{jlmn},
\end{equation}
and
\begin{equation}\label{nabla-psi}
  \nabla_m\psi_{ijkl}=-\left(T_{mi}\varphi_{jkl}-T_{mj}\varphi_{ikl}-T_{mk}\varphi_{jil}-T_{ml}\varphi_{jki}\right).
\end{equation}
Then the full torsion tensor $T_{ij}$ is related with the torsion forms by the following (see \cite{Kar}):
\begin{equation*}
 T_{ij}=\frac{\tau_0}4g_{ij}-(\tau_1^{\#}\lrcorner\varphi)_{ij}-(\bar{\tau}_3)_{ij}-\frac 12(\tau_2)_{ij},
\end{equation*}
where $(\tau_1^{\#}\lrcorner\varphi)_{ij}=(\tau_1)^l\varphi_{lij}$ and $(\bar{\tau}_3)_{ij}$ is the traceless symmetric 2-tensor such that $\tau_3=i_{\varphi}(\bar{\tau}_3)$.

If $\varphi$ is closed, then \eqref{dvarphi} implies that $\tau_0,\tau_{1}$ and $\tau_3$ are all zero, so the only non-zero torsion form is $\tau_2=\frac 12(\tau_2)_{ij}dx^i\wedge dx^j$. The full torsion tensor satisfies $T_{ij}=-T_{ji}=-\frac 12(\tau_2)_{ij}$, which is a skew-symmetric $2$-tensor. In the majority of this article, we will only consider closed $\GG_2$ structures and will thus write $\tau=\tau_2$ for simplicity.

To study the Laplacian flow we have calculated $\Delta_{\varphi}\varphi$ for closed $\varphi$ in \cite{Lotay-Wei}:
 \begin{equation}
   \Delta_{\varphi}\varphi=d\tau=i_{\varphi}(h)
 \end{equation}
for a symmetric $2$-tensor $h$ satisfying
\begin{equation}\label{hodge-Lap-varp-3}
h_{ij}=-(\nabla_mT_{ni})\varphi_j^{\,\,\,mn}-\frac 13|T|^2g_{ij}-T_{i}^{\,\,l}T_{lj}.
\end{equation}
Thus the Laplacian flow \eqref{Lap-flow-def} can be rewritten as
\begin{equation}\label{flow-general}
  \frac{\pt}{\pt t}\varphi=i_{\varphi}(h).
\end{equation}

Under the flow \eqref{flow-general}, the associated metric $g(t)$ of $\varphi(t)$ evolves by (see \cite{Kar})
\begin{equation*}
  \frac{\pt}{\pt t}g(t)=2h(t).
\end{equation*}
Then substituting \eqref{hodge-Lap-varp-3}, we have
\begin{equation}\label{flow-g1}
  \frac{\pt}{\pt t}g_{ij}=-2(\nabla_kT_{li})\varphi_{j}^{\,\,kl}-\frac 23|T|^2g_{ij}-2T_i^{\,\,k}T_{kj}.
\end{equation}
Since $d\varphi=0$, the Ricci tensor of the associated metric $g$ is equal to (see \cite{Kar})
\begin{equation}\label{ric}
  Ric_{ij}=(\nabla_kT_{li})\varphi_{j}^{\,\,kl}-T_i^{\,\,k}T_{kj},
\end{equation}
 and so we can recover the following formula from \cite{bryant2005}:
\begin{equation}\label{flow-g2}
  \frac{\pt}{\pt t}g_{ij}=-2Ric_{ij}-\frac 23|T|^2g_{ij}-4T_i^{\,\,k}T_{kj}.
\end{equation}
Thus the leading term of the metric flow corresponds to the Ricci flow. 
 From \eqref{flow-g1}, the volume form $vol_g$ evolves by
\begin{equation}\label{evl-volumeform}
  \frac{\pt}{\pt t}vol_g=\frac 12tr_g(\frac{\pt}{\pt t}g(t))vol_g=\frac 23|T|^2vol_g.
\end{equation}
Therefore along the Laplacian flow for closed $\GG_2$ structures, the volume of $M$ with respect to $g$ will increase
 (in fact, the volume form increases pointwise).

 Notice, as in \cite{bryant2005}, that taking the trace of \eqref{ric} gives
\begin{equation}\label{eq-scal}
R=-|T|^2,
\end{equation}
 where $R$ is the scalar curvature.  So, 7-manifolds with closed $\GG_2$ structures necessarily have metrics with non-positive scalar curvature, but this is a weak curvature condition amenable to the h-principle.

In \cite{Lotay-Wei}, we also calculated the evolution of the full torsion tensor $T$,
\begin{equation}\label{evl-torsion-1}
  \frac{\pt}{\pt t}T=\Delta T+Rm*T+Rm*T*\psi+\nabla T*T*\varphi+T*T*T,
\end{equation}
where we will continually use $*$ to mean some contraction using the metric $g(t)$ associated with $\varphi(t)$, and the evolution equation of the  Riemann curvature tensor as
\begin{equation}\label{evl-Rm}
  \frac{\pt}{\pt t} Rm=\Delta Rm+Rm*Rm+Rm*T*T+\nabla^2T*T+\nabla T*\nabla T.
\end{equation}
Here $\Delta$ denotes the usual (analyst's) Laplacian, which is a negative operator.  Notice that $\Delta_{\varphi}$ on the other hand is the Hodge Laplacian and so is a positive operator.

We end this section with the following useful estimate on the existence time of the Laplacian flow \eqref{Lap-flow-def} from \cite{Lotay-Wei}.
\begin{prop}
\label{prop-prelim}
If $\varphi_0$ is a closed $\GG_2$ structure on a compact 7-manifold $M$ and $K\in\R$ so that
\begin{equation}\label{eq-Lambda}
  \Lambda(x)=\left(|Rm(x)|^2+|\nabla T(x)|^2\right)^{\frac 12}\leq K,
\end{equation}
for all $x\in M$, then the unique solution $\varphi(t)$ of the Laplacian flow \eqref{Lap-flow-def} will exist at least for time $t\in [0,c/{K}]$, where $c>0$ is a uniform  constant, independent of $\varphi_0$.
\end{prop}

\noindent Observe that a
 bound on $Rm$ gives a bound on $T$ for closed $\GG_2$ structures by \eqref{eq-scal}.

\section{Estimates for \texorpdfstring{$t\leq 1$}{t<1}}\label{sec-estimate}
In this section, we show that if $\bar{\varphi}$ is a torsion-free $\GG_2$ structure and we are given a neighbourhood $U$ of $\bar{\varphi}$, then we can find a smaller neighbourhood $U_0$ of $\bar{\varphi}$ so that  if the initial value $\varphi_0\in [\bar{\varphi}]_+$  lies in $U_0$, then the Laplacian--DeTurck flow \eqref{deturck-flow}  will exist at least for $t\in [0,1]$ and remain in $U$. 

\begin{lem}[Estimate for $t\leq 1$]\label{lem-est-t-leq1}
Let $(M,\bar{\varphi})$ be a compact $\GG_2$ manifold. 
For any $k\geq 3$, $\epsilon>0$, there exists a constant $\delta=\delta(M,\bar{\varphi},k,\epsilon)>0$ such that: if $\varphi_0\in [\bar{\varphi}]_+$ and $\|\varphi_0-\bar{\varphi}\|_{C^{k+4}_{\bar{g}}}<\delta$, then the solution $\tilde{\varphi}(t)$ to the Laplacian--DeTurck flow \eqref{deturck-flow} 
exists on $[0,1]$ and satisfies
\begin{equation*}
 \|\tilde{\varphi}(t)-\bar{\varphi}\|_{C^{k}_{\bar{g}}}<\epsilon \quad \forall t\in [0,1],
\end{equation*}
where the $C^k_{\bar{g}}$-norms are defined via the metric $\bar{g}$ determined by $\bar{\varphi}$.
\end{lem}
\proof
Let $k\geq 3$ and $\epsilon>0$.   By making $\epsilon>0$ smaller, we can assume that if $\varphi$ is a $\GG_2$ structure satisfying $\|\varphi-\bar{\varphi}\|_{C^k_{\bar{g}}}\leq \epsilon$ and $g$ is the associated metric of $\varphi$, then the $C^k$-norms via $g$ and  $\bar{g}$ differ at most by a factor $2$.

Recall  $\Lambda$ from \eqref{eq-Lambda}. Define a finite constant $K=K(\epsilon,k)\geq 0$ by:
\begin{align}
  K= &\sup\{ |\nabla^iRm(x)|^2 +|\nabla^{i+1}T(x)|^2\,:\,\|\varphi-\bar{\varphi}\|_{C^{k+4}_{\bar{g}}}\leq \epsilon, x\in M, i\leq k+2\}\nonumber\\
  & +\sup\{|\Lambda(x)|\,:\, \|\varphi-\bar{\varphi}\|_{C^k_{\bar{g}}}\leq \epsilon, x\in M\},\label{eq-K}
\end{align}
where $\nabla$, $Rm$ and $T$ denote the Levi-Civita connection, Riemann curvature tensor and torsion tensor of $\varphi$.

Let $\varphi_0$ satisfy $\|\varphi_0-\bar{\varphi}\|_{C^{k+4}_{\bar{g}}}\leq\epsilon$
and let $\varphi(t)$ be the solution to the Laplacian flow \eqref{Lap-flow-def} starting at $\varphi_0$.
Suppose that $\varphi(t)$ exists on some time interval $[0,\eta]$ (where $\eta>0$ by Theorem \ref{thm-bryant-xu}). 
From \eqref{evl-torsion-1}-\eqref{evl-Rm} and calculations in \cite[\S 4]{Lotay-Wei}, we can compute the evolution of the pointwise norms $|\nabla^i Rm|$ and
$|\nabla^{i+1}T|$ with respect to $g(t)$.

First we have
\begin{align}\label{dt-na^k-Lamb}
  &\frac{\pt}{\pt t}(|\nabla^iRm|^2 +|\nabla^{i+1}T|^2)\nonumber\\&
  \leq  \Delta(|\nabla^iRm|^2 +|\nabla^{i+1} T|^2)-2|\nabla^{i+1}Rm|^2-2|\nabla^{i+2}T|^2\nonumber\\
  &+\sum_{j=0}^i\nabla^iRm*\nabla^{i-j}Rm*\nabla^j(Rm+T*T)\nonumber\\
  &+\sum_{j=0}^{i+1}\nabla^iRm*\nabla^jT*\nabla^{i+2-j}T+\sum_{j=0}^{i+1}\nabla^{i+1} T*\nabla^{i+1-j}T*\nabla^j(Rm+T*T)\nonumber\\
  &+\sum_{j=0}^{i+1}\nabla^{i+1} T*(\nabla^{i+1-j}(Rm*T)*\nabla^j\psi+\nabla^{i+1-j}(\nabla T*T)*\nabla^j\varphi).
\end{align}
Note that the pointwise norms of $\nabla^j\varphi$ and $\nabla^j\psi$ can be estimated by a constant depending on the upper bounds of $|T|,|\nabla T|,\cdots, |\nabla^{j-1}T|$ 
and that \eqref{dt-na^k-Lamb} for $i=0$ gives an evolution inequality for $\Lambda^2$.

By definition of $K$ in \eqref{eq-K} we have that if the Laplacian flow $\varphi(t)$  starting from $\varphi_0$ exists on the time interval $[0,\eta]$ and $\|\varphi(t)-\bar{\varphi}\|_{C^k_{\bar{g}}}\leq \epsilon$ for $\forall t\in [0,\eta]$, then
\begin{equation}\label{K-def}
  \left\{\begin{array}{l}
                           \sup_{x\in M}\Lambda(x,t)\leq K,\quad \forall t\in [0,\eta],\\
                           ~\\
                           |\nabla^iRm(x,0)|^2 +|\nabla^{i+1}T(x,0)|^2\leq K,\forall x\in M, i\leq k+2.
                         \end{array}
 \right.
\end{equation}
Applying the maximum principle to \eqref{dt-na^k-Lamb} and using
 an  induction argument, we claim that if $\eta\leq 1$, then there is a constant $\tilde{K}=\tilde{K}(K,k)<\infty$  such that
 \begin{equation}\label{na-Lamb}
   |\nabla^iRm(x,t)|^2 +|\nabla^{i+1}T(x,t)|^2\leq \tilde{K},\forall x\in M, t\in [0,\eta], i\leq k+2.
 \end{equation}
We now prove this claim. 

For $i=1$, using \eqref{eq-scal}  we can derive from \eqref{dt-na^k-Lamb} that
\begin{align}\label{evl-d-Rm-T-0}
  \frac{\pt}{\pt t}(|\nabla Rm|^2 +|\nabla^2T|^2) &\leq \Delta (|\nabla Rm|^2 +|\nabla^2T|^2)-2|\nabla^2Rm|^2-2|\nabla^3T|^2\nonumber\\
  &+C|\nabla^3 T||Rm|^{\frac{1}{2}}(|\nabla Rm|+|\nabla^2T|)\nonumber\\
  &+C|\nabla^2 T||\nabla^2 Rm||Rm|^{\frac{1}{2}}\nonumber\\
   &+C(|\nabla Rm|^2+|\nabla^2T|^2)(|\nabla T|+|Rm|)\nonumber\\
  &+C(|\nabla Rm|+|\nabla^2T|)|Rm|^{\frac 12}(|Rm|^2+|\nabla T|^2),
\end{align}
where $C$ are uniform constants.  Young's inequality implies that $2ab\leq \frac{1}{\epsilon}a^2+\epsilon b^2$ for all $\epsilon>0$ and $a,b\geq 0$, so for all $\epsilon>0$ we have
\begin{align*}
2|\nabla^3 T||Rm|^{\frac{1}{2}}(|\nabla Rm|+|\nabla^2T|)&\leq 
\frac{2}{\epsilon}|Rm|(|\nabla Rm|^2+|\nabla^2 T|^2)+\epsilon|\nabla^3 T|^2,\\
2|\nabla^2 T||\nabla^2 Rm||Rm|^{\frac{1}{2}}&\leq \frac{1}{\epsilon}|Rm||\nabla^2 T|^2+\epsilon |\nabla^2 Rm|^2.
\end{align*}
Substituting these estimates into \eqref{evl-d-Rm-T-0}  and choosing $\epsilon$ sufficiently small then yields the parabolic inequality
\begin{align}\label{evl-d-Rm-T}
  \frac{\pt}{\pt t}(|\nabla Rm|^2 +|\nabla^2T|^2) &\leq \Delta (|\nabla Rm|^2 +|\nabla^2T|^2)\nonumber\\
   &+C(|\nabla Rm|^2+|\nabla^2T|^2)(|\nabla T|+|Rm|)\nonumber\\
  &+C(|\nabla Rm|+|\nabla^2T|)|Rm|^{\frac 12}(|Rm|^2+|\nabla T|^2).
\end{align}
 If we denote $$\Lambda_1(t)=\max_{x\in M}(|\nabla Rm(x,t)|^2 +|\nabla^2T(x,t)|^2),$$ then applying the maximum principle to \eqref{evl-d-Rm-T} and using the first line of the condition \eqref{K-def}, we have
\begin{align}\label{evl-Lambda_1}
  \frac d{dt}\Lambda_1(t)\leq  & C_1K\Lambda_1(t)+C_2K^2,
\end{align}
where $C_1,C_2$ are two uniform constants. Then integrating \eqref{evl-Lambda_1} yields that
\begin{equation*}
\Lambda_1(t)+\frac {C_2}{C_1}K\leq e^{C_1Kt}(\Lambda_1(0)+\frac {C_2}{C_1}K).
\end{equation*}
Since $t\leq 1$ and $\Lambda_1(0)\leq K$ by the second line of \eqref{K-def}, the above inequality implies that
\begin{equation}\label{na-Lamb-1}
\Lambda_1(t)\leq e^{C_1K}K+(e^{C_1K}-1)\frac {C_2}{C_1}K:=\tilde{K}.
\end{equation}
This proves the $i=1$ case of \eqref{na-Lamb}. 

For $1<i\leq k+2$, the conclusion \eqref{na-Lamb} follows from an induction argument starting from \eqref{na-Lamb-1}.  More concretely, suppose that there is some $1<l\leq k+2$ so that \eqref{na-Lamb}  
  holds for all $1\leq i<l$.  Consider \eqref{dt-na^k-Lamb} for $i=l$.  By the inductive hypothesis, any term on the right-hand side of \eqref{dt-na^k-Lamb} which can be estimated by $|\nabla^iRm|$ and $|\nabla^{i+1}T|$ for $i<l$ can be bounded by a constant depending only on $\tilde{K}$.  With this observation, \eqref{dt-na^k-Lamb} for $i=l$ then yields
\begin{align*}
\frac{\pt}{\pt t}(|\nabla^lRm|^2 &+|\nabla^{l+1}T|^2)\\
&  \leq  \Delta(|\nabla^lRm|^2 +|\nabla^{l+1} T|^2)-2|\nabla^{l+1}Rm|^2-2|\nabla^{l+2}T|^2\\
  &+C|\nabla^{l+2}T|(|\nabla^lRm|+|\nabla^{l+1}T|)+C|\nabla^{l+1}T||\nabla^{l+1}Rm|\\
  &+C(|\nabla^lRm|^2+|\nabla^{l+1}T|^2)+C(|\nabla^lRm|+|\nabla^{l+1}T|),
\end{align*}
where the constants $C$ depend on $\tilde{K}$.  
Using Young's inequality as before we obtain the parabolic inequality
\begin{align*}
\frac{\pt}{\pt t}(|\nabla^lRm|^2 &+|\nabla^{l+1}T|^2)
 \leq  \Delta(|\nabla^lRm|^2 +|\nabla^{l+1} T|^2)\\
  &\qquad\qquad+C(|\nabla^lRm|^2+|\nabla^{l+1}T|^2)+C(|\nabla^lRm|+|\nabla^{l+1}T|),
\end{align*}
to which we may apply the maximum principle to conclude that \eqref{na-Lamb} holds for $i=l$ as required.

Next, using calculations in \cite[\S 4]{Lotay-Wei} again, we have
\begin{align}\label{dt-na^k-T}
  &\frac{\pt}{\pt t}|\nabla^{i}T|^2\leq  \Delta|\nabla^{i} T|^2-2|\nabla^{i+1}T|^2 
  +\sum_{j=0}^{i}\nabla^{i} T*\nabla^{i-j}T*\nabla^j(Rm+T*T)\nonumber\\
  &+\sum_{j=0}^{i}\nabla^{i} T*\nabla^{i-j}(Rm*T)*\nabla^j\psi 
+\sum_{j=0}^{i}\nabla^{i} T*\nabla^{i-j}(\nabla T*T)*\nabla^j\varphi.
\end{align}
We now choose $\tilde{\epsilon}>0$, which we will impose smallness conditions on later.  Applying the maximum principle to \eqref{dt-na^k-T} and using \eqref{na-Lamb}, there is a constant $\tilde{\delta}=\tilde{\delta}(K,\tilde{\epsilon},k)>0$ such that
\begin{align}\label{na-i-T}
 |\nabla^{i}T(x,0)|\leq \tilde{\delta}\quad &\forall x\in M, i\leq k+2
  \nonumber\\
  &
  \Rightarrow |\nabla^{i}T(x,t)|\leq \tilde{\epsilon}\quad \forall x\in M, t\in [0,\eta], i\leq k+2.
\end{align}
Note that in \eqref{na-i-T} we have control of one less derivative of $T$ than in \eqref{na-Lamb}, because on the right-hand side of \eqref{dt-na^k-T} there are terms involving $\nabla^iRm$ and $\nabla^{i+1}T$ which need to be estimated.

Since $\|\varphi_0-\bar{\varphi}\|_{C^{k+4}_{\bar{g}}}\leq \epsilon$, by \eqref{eq-K} we have $\Lambda(x,0)\leq K<\infty$,  and hence we see from Proposition \ref{prop-prelim} that the maximal existence time $T_0$ of $\varphi(t)$ has a uniform lower bound $T_0\geq c/{K}>0$ depending only on $K$.  We know that the solution $\tilde{\varphi}(t)$ to the Laplacian--DeTurck flow \eqref{deturck-flow}, with the same initial value $\varphi_0$,
is given by
\[
\tilde{\varphi}(t)=\phi(t)^*\varphi(t)
\]
for a family of diffeomorphisms $\phi(t)$ on $M$ generated by vector fields $V(\varphi(t))$ determined by $\varphi(t)$.  Hence $\tilde{\varphi}(t)$ will also exist at least on $[0,T_0)$.

 We also know that $V(\varphi(t))$ depends on $\varphi(t)$ and its first derivative, 
 so a bound on $|\nabla^i\varphi(t)|_{g(t)}$ for $i\leq k+2$ will yield a bound on $|\tilde{\nabla}^i\tilde{\varphi}|_{\tilde{g}(t)}$
  for $i\leq k$, where $\tilde{g}(t)$ and $\tilde{\nabla}$ are the metric and Levi-Civita connection determined by $\tilde{\varphi}(t)$.  
  Therefore, since $\varphi(t)$ is bounded in $C^{k+4}_{g}$ we have that $\tilde{\varphi}(t)$ is bounded in $C^{k+2}_{\tilde{g}}$ and thus 
\begin{equation*}
  |\tilde{\nabla}^i\widetilde{Rm}|_{\tilde{g}}^2+|\tilde{\nabla}^{i+1}\tilde{T}|_{\tilde{g}}^2 
\end{equation*}
is finite for $i\leq k $, where $\widetilde{Rm}$ and $\tilde{T}$ are the Riemann curvature and torsion tensors of $\tilde{\varphi}(t)$ respectively. By diffeomorphism invariance, \eqref{na-Lamb} and \eqref{na-i-T} show that
\begin{align}
  |\tilde{\nabla}^i\widetilde{Rm}|_{\tilde{g}}^2+|\tilde{\nabla}^{i+1}\tilde{T}|_{\tilde{g}}^2&=|\nabla^iRm|_g^2+|\nabla^{i+1}T|_g^2\nonumber\\
  &\leq \tilde{K},\quad \forall x\in M, t\in[0,\eta],i\leq k ,\label{tna-Lambda}\\
  |\tilde{\nabla}^{i}\tilde{T}|_{\tilde{g}}&=|\nabla^{i}T|_g\nonumber\\
  &\leq\tilde{\epsilon}, \quad\;\,\, \forall x\in M, t\in[0,\eta], i\leq k+1 \label{tna-i-T}.
\end{align}


Now let $\delta\in(0,\epsilon)$ be sufficiently small so that for all $x\in M$, $i\leq k+2$ we have
\begin{align*}
    \|\varphi-\bar{\varphi}\|_{C^{k+4}_{\bar{g}}}\leq \delta  \Rightarrow
    \|\varphi-\bar{\varphi}\|_{C^k_{\bar{g}}}\leq \epsilon/4
    \textrm{ and }|\nabla^{i}T(x)|\leq \tilde{\delta}.
\end{align*}
Assume that $\|\varphi_0-\bar{\varphi}\|_{C^{k+4}_{\bar{g}}}\leq \delta<\epsilon$, let $\eta\leq T_0$ be the maximal existence time of the  Laplacian--DeTurck flow \eqref{deturck-flow} with the initial value $\varphi_0$ such that
\begin{equation*}
  \|\tilde{\varphi}(t)-\bar{\varphi}\|_{C^k_{\bar{g}}}< \epsilon \quad \forall t\in [0,\eta).
\end{equation*}
If $\eta<1$, then
\begin{align}\label{est-shorttime}
 \|\tilde{\varphi}(t)-\bar{\varphi}\|_{C^k_{\bar{g}}}\leq & \|\varphi_0-\bar{\varphi}\|_{C^k_{\bar{g}}}+\int_0^{t}\frac d{ds}\|\tilde{\varphi}(s)-\bar{\varphi}\|_{C^k_{\bar{g}}}ds \nonumber\\
  \leq  & \frac{\epsilon}4+\int_0^{\eta}\|\Delta_{\tilde{\varphi}(s)}\tilde{\varphi}(s)+\mathcal{L}_{V(\tilde{\varphi}(s))}\tilde{\varphi}(s)\|_{C^k_{\bar{g}}}ds\nonumber\\
  \leq  & \frac{\epsilon}4+2\int_0^{\eta}\|\Delta_{\tilde{\varphi}(s)}\tilde{\varphi}(s)+\mathcal{L}_{V(\tilde{\varphi}(s))}\tilde{\varphi}(s)\|_{C^k_{\tilde{g}(s)}}ds\nonumber\displaybreak[0]\\
  \leq  & \frac {\epsilon}4+C\int_0^{\eta}(\|\tilde{\nabla} \tilde{T}(x,s)\|_{C^k_{\tilde{g}(s)}}+\|\tilde{T}(x,s)\|^2_{C^k_{\tilde{g}(s)}})\nonumber\\
  \leq &\frac {\epsilon}4+2C\tilde{\epsilon}
\end{align}
for all $t\in[0,\eta]$, where in the final step we used $\eta<1$.  Here 
 in the fourth inequality of \eqref{est-shorttime} we used the facts $\Delta_{\tilde{\varphi}}\tilde{\varphi}=d\tilde{\tau}$ and the vector field $V(\tilde{\varphi})$ involves the first order covariant derivative of $\tilde{\varphi}$, and in the last inequality we used \eqref{tna-i-T}, noting that we need a small bound on $\tilde{T}$ in $C^{k+1}_{\tilde{g}}$.  
 Since $C$ is a uniform constant, by choosing $\tilde{\epsilon}$ sufficiently small we can ensure that $2C\tilde{\epsilon}<\frac{\epsilon}{4}$, and hence
  $\|\tilde{\varphi}(t)-\bar{\varphi}\|_{C^k_{\bar{g}}}\leq\frac{\epsilon}{2}$,  contradicting the maximality of $\eta$.

Therefore we have $\eta\geq 1$ as desired. 
\endproof

\section{Estimates for \texorpdfstring{$t\geq 1$}{t>1}}\label{sec:t>1}

In this section, we now turn to studying the long-time estimates for solutions $\tilde{\varphi}(t)$ to the Laplacian--DeTurck flow \eqref{deturck-flow}; i.e.~for $t\geq 1$.  We will show that
the flow exponentially decays to the torsion-free $\GG_2$ structure
in $L^2$.  We then obtain further integral estimates which
allow us to deduce, through Sobolev embedding, that the $C^k_{\bar{g}}$-norm of $\tilde{\varphi}(t)-\bar{\varphi}$
is uniformly controlled by the $C^{k+4}_{\bar{g}}$-norm of $\varphi_0-\bar{\varphi}$.

We suppose we are in
the situation of Lemma \ref{lem-est-t-leq1} and use the notation introduced there: in particular, we fix $k\geq 3$ and $\epsilon>0$ and obtain $\delta>0$ so that the flow exists at least up to time $1$ in
$[\bar{\varphi}]_+$ where $\bar{\varphi}$ is torsion-free.

Since $\tilde{\varphi}(t)$ lies in the cohomology class $[\bar{\varphi}]$, we can write $\tilde{\varphi}(t)=\bar{\varphi}+\theta(t)$ with $\theta(t)$ taking values in exact forms.  We are assuming that $\theta_0=\theta(0)$ satisfies $\|\theta_0\|_{C^{k+4}_{\bar{g}}}<\delta$.  Let $\eta$ be the maximum time such that the solution $\tilde{\varphi}(t)$ of $\eqref{deturck-flow}$ starting at $\varphi_0$ exists and $\|\tilde{\varphi}(t)-\bar{\varphi}\|_{C^k_{\bar{g}}}<\epsilon$. Clearly, from Lemma \ref{lem-est-t-leq1} we know that $\eta\geq 1$.

Since we shall focus on the Laplacian--DeTurck flow \eqref{deturck-flow} in this section, we now recall the definition of the vector field $V(\varphi)$ which appears in the flow.  Let $S=\nabla-\bar{\nabla}$, which is a tensor taking values in $TM\otimes S^2T^*M$.  Let $A$ be the constant in \cite[Lemma 1.7]{bryant-xu2011}.  (Although $A$ is not calculated in \cite{bryant-xu2011}, it arises purely from $\GG_2$ representation theory and so by working on $\R^7$ one can find that $A=-\frac{1}{8}$ in the conventions of \cite{bryant-xu2011}.) As in\cite{bryant-xu2011}, we let $V(\varphi)$ be
\begin{equation}\label{vector-V}
  V(\varphi)=-5V_1(\varphi)-V_2(\varphi),
\end{equation}
with $V_1,V_2$ given locally by
\begin{equation*}
  V_1(\varphi)=\frac 17g^{pq}S_{pq}^ie_i,\qquad V_2(\varphi)=2A(g^{kj}S_{ik}^ie_j+5V_1(\varphi)).
\end{equation*}
 In particular, observe that $V(\bar{\varphi})=0$. We let $\mathcal{A}_{\bar{\varphi}}(\varphi)=\Delta_{\varphi}\varphi+\mathcal{L}_{V(\varphi)}\varphi$. It was proved in \cite{bryant-xu2011} that the linearization of $\mathcal{A}_{\bar{\varphi}}(\varphi)$ at $\varphi$ applied to a closed variation $\theta=\pt_s|_{s=0}\varphi_s$ (so $d\theta=0$) is
\begin{equation}\label{linearize-A}
  D_{\varphi}\mathcal{A}_{\bar{\varphi}}(\varphi)
  \theta=-\Delta_{\varphi}\theta+d\Phi(\theta),
\end{equation}
where $\Phi(\theta)$ is an algebraic linear operator on $\theta$ with coefficients depending on the torsion of $\varphi$.  Notice in \eqref{linearize-A} the surprising fact that the sign of the Laplacian has changed in the linearisation, which indicates the parabolicity of the flow discussed in the introduction. We should also note that \eqref{linearize-A} was derived only assuming $\bar{\varphi}$ was a background $\GG_2$ structure and not necessarily a torsion-free one.

We now show that we can approximate the evolution of $\theta(t)=\tilde{\varphi}(t)-\bar{\varphi}$.
\begin{lem}\label{lem-linearize}
There exists $\bar{\epsilon}=\bar{\epsilon}(\bar{\varphi})>0$ such that if $\epsilon\in(0,\bar{\epsilon})$, then for all time $t$ for which the solution $\tilde{\varphi}(t)$ of $\eqref{deturck-flow}$ exists and $\|\theta(t)\|_{C^k_{\bar{g}}}<\epsilon$, we have
\begin{align}\label{flow-theta}
  \frac{\pt}{\pt t}\theta(t) =&-\Delta_{\bar{\varphi}}\theta+dF(\bar{\varphi},\tilde{\varphi}(t),\theta(t),\bar{\nabla}\theta(t)),
\end{align}
where $F$ is a $2$-form which is smooth in the first two arguments and linear in the last two arguments.
\end{lem}
\proof
Recall that $\tilde{\varphi}(t)$ is a solution to the Laplacian--DeTurck flow \eqref{deturck-flow} and $\bar{\varphi}$ is a stationary solution of \eqref{deturck-flow}, so $\mathcal{A}_{\bar{\varphi}}(\bar{\varphi})=0$.  For all $t$ for which $\tilde{\varphi}(t)$ exists and $\|\theta(t)\|_{C^k_{\bar{g}}}=\|\tilde{\varphi}(t)-\bar{\varphi}\|_{C^k_{\bar{g}}}<\epsilon$, we have
\begin{align}
  \frac{\pt}{\pt t}\theta(t) =& \frac{\pt }{\pt t}(\tilde{\varphi}(t)-\bar{\varphi}) =\mathcal{A}_{\bar{\varphi}}(\tilde{\varphi}(t))-\mathcal{A}_{\bar{\varphi}}(\bar{\varphi})\nonumber\\
  =&D_{\bar{\varphi}}\mathcal{A}_{\bar{\varphi}}(\tilde{\varphi})\theta+\tilde{F}=-\Delta_{\bar{\varphi}}\theta+d\Phi(\theta)+\tilde{F}\label{dt-theta}
\end{align}
where $\tilde{F}$ is the remaining higher order term. We choose $\epsilon$ small enough such that $\tilde{F}$ is small compared to the linear part. By studying the details of the proof of \eqref{linearize-A} in \cite{bryant-xu2011}, we can conclude that if we linearize $\mathcal{A}_{\bar{\varphi}}(\tilde{\varphi})$ at the torsion-free $\GG_2$ structure $\bar{\varphi}$, then $d\Phi(\theta)=0$, so
\begin{equation}\label{lin-Delta}
D_{\bar{\varphi}}\mathcal{A}_{\bar{\varphi}}(\tilde{\varphi})\theta=-\Delta_{\bar{\varphi}}\theta.
\end{equation}

We next estimate the remaining term $\tilde{F}$. Since $\theta(t)$ is small in the $C^k_{\bar{g}}$-norm, we can write 
\begin{equation}\label{theta-decomp}
\theta=3f^0\bar{\varphi}+*_{\bar{\varphi}}(f^1\wedge\bar{\varphi})+f^3
\end{equation} for $f^0\in \Omega^0(M)$, $f^1\in \Omega^1(M)$ and $f^3\in \Omega^3_{27}(M)$ (with respect to $\bar{\varphi}$) with small $C^k_{\bar{g}}$-norm. Then by  \cite{joyce2000} we have
\begin{align}
  *_{\tilde{\varphi}}\tilde{\varphi}= & *_{\bar{\varphi}}\bar{\varphi}+4f^0*_{\bar{\varphi}}\bar{\varphi}+f^1\wedge\bar{\varphi}-*_{\bar{\varphi}}f^3+Q(\bar{\varphi},\tilde{\varphi},\theta),
\end{align}
where $Q$ is smooth in its first two arguments, and the principal part of $Q$ is quadratic in the third argument. As $\tilde{\varphi}(t)$ is closed,   we have
\begin{align*}
  \Delta_{\tilde{\varphi}}\tilde{\varphi}= *_{\tilde{\varphi}}d*_{\tilde{\varphi}}d\tilde{\varphi}-d*_{\tilde{\varphi}}d*_{\tilde{\varphi}}\tilde{\varphi}
  = -d*_{\tilde{\varphi}}d*_{\tilde{\varphi}}\tilde{\varphi}.
\end{align*}
We linearize $\Delta_{\tilde{\varphi}}\tilde{\varphi}$ at $\bar{\varphi}$,
\begin{align*}
  D_{\bar{\varphi}}\Delta_{\tilde{\varphi}}\tilde{\varphi}(\theta)  =& -d(\dot{*}_{\tilde{\varphi}})d*_{\bar{\varphi}}\bar{\varphi} -d*_{\bar{\varphi}}d(4f^0*_{\bar{\varphi}}\bar{\varphi}+f^1\wedge\bar{\varphi}-*_{\bar{\varphi}}f^3)\\
  =& -d*_{\bar{\varphi}}d(4f^0*_{\bar{\varphi}}\bar{\varphi}+f^1\wedge\bar{\varphi}-*_{\bar{\varphi}}f^3),
\end{align*}
where we used the fact that $\bar{\varphi}$ is coclosed (as it is torsion-free) to eliminate the first term involving the linearisation of $*_{\tilde{\varphi}}$.
We then have
\begin{align}
  \Delta_{\tilde{\varphi}}\tilde{\varphi}-\Delta_{\bar{\varphi}}\bar{\varphi}&-D_{\bar{\varphi}}\Delta_{\tilde{\varphi}}\tilde{\varphi}(\theta) \nonumber\\ &=-d*_{\tilde{\varphi}}d*_{\tilde{\varphi}}\tilde{\varphi}+d*_{\bar{\varphi}}d(4f^0*_{\bar{\varphi}}\bar{\varphi}+f^1\wedge\bar{\varphi}-*_{\bar{\varphi}}f^3)\nonumber\\
  &=  -d(*_{\tilde{\varphi}}-*_{\bar{\varphi}})d*_{\tilde{\varphi}}\tilde{\varphi}-d*_{\bar{\varphi}}dQ(\bar{\varphi},\tilde{\varphi},\theta)\nonumber\\
  &=dF_1(\bar{\varphi},\tilde{\varphi},\theta,\bar{\nabla}\theta),\label{pf-linerar-1}
\end{align}
where we denote $$F_1(\bar{\varphi},\tilde{\varphi},\theta,\bar{\nabla}\theta)=-(*_{\tilde{\varphi}}-*_{\bar{\varphi}})d\!*_{\tilde{\varphi}}\!\tilde{\varphi}-*_{\bar{\varphi}}dQ(\bar{\varphi},\tilde{\varphi},\theta).$$ Notice that  $F_1$ is dominated by $\theta\bar{\nabla}\theta$ for $\epsilon$ small.

Let $\tilde{g}=\bar{g}+h$ be the metric determined by $\tilde{\varphi}$. To calculate the difference between $\mathcal{L}_{V(\tilde{\varphi})}\tilde{\varphi}$ and its linearization $D_{\bar{\varphi}}(\mathcal{L}_{V(\tilde{\varphi})}\tilde{\varphi})(\theta)$, we first calculate the difference between $V_1(\tilde{\varphi})=\tilde{g}^{pq}S_{pq}^ie_i$ and its linearization. We can express $V_1(\tilde{\varphi})$ in local coordinates as follows:
\begin{align}
  V_1(\tilde{\varphi}) =& (\bar{g}+h)^{pq}(\Gamma_{pq}^i(\bar{g}+h)-\Gamma_{pq}^i(\bar{g}))e_i \nonumber\\
  = & \frac 12(\bar{g}+h)^{pq}\biggl((\bar{g}+h)^{ij}\left(\pt_p(\bar{g}+h)_{jq}+\pt_q(\bar{g}+h)_{jp}-\pt_j(\bar{g}+h)_{pq}\right)\nonumber\\
  &\quad\qquad\qquad-\bar{g}^{ij}(\pt_p(\bar{g})_{jq}+\pt_q(\bar{g})_{jp}-\pt_j(\bar{g})_{pq})\biggr)e_i.\label{V1-phi-eq}
\end{align}
Note that $V_1(\tilde{\varphi})$ is a tensor and its linearization is:
\begin{equation}\label{lin-V}
  V_{1*}(\theta)=\frac 12\bar{g}^{pq}\bar{g}^{ij}\left(\bar{\nabla}_ph_{jq}+\bar{\nabla}_qh_{jp}-\bar{\nabla}_jh_{pq}\right)e_i.
\end{equation}
It is well-known (c.f.~\cite{bryant2005,joyce2000}) that the metric $\tilde{g}$ of $\tilde{\varphi}$ varies infinitesimally by 
\begin{equation}\label{h-decomp}
h=2f^0\bar{g}+\frac 12j_{\bar{\varphi}}(f^3),
\end{equation} where $j_{\bar{\varphi}}$ is the linear map defined in \eqref{j-varphi-def}.  Equations \eqref{theta-decomp} and \eqref{h-decomp} show that there is a linear relation between $h$ and $\theta$.  Therefore,  using the fact that $V_1(\bar{\varphi})=0$, we may write
\begin{align*}
  V_1(\tilde{\varphi})-V_1(\bar{\varphi})-V_{1*}(\theta)=& F_2(\bar{g},\tilde{g},h,\bar{\nabla} h)=F_2(\bar{\varphi},\tilde{\varphi},\theta,\bar{\nabla}\theta),
\end{align*}
where $F_2$ is linear in $\theta$ and $\bar{\nabla}\theta$.  Moreover, it follows from \eqref{V1-phi-eq} and \eqref{lin-V}, and the relation between $h$ and $\theta$, that $F_2$  can be expressed in the form $f(\bar{\varphi},\tilde{\varphi})\theta\bar{\nabla}\theta$, where $f$ is a tensor depending on $\bar{\varphi}$ and $\tilde{\varphi}$.

By the same method, we can express the difference between the vector field $\tilde{g}^{kj}S_{ik}^ie_j$ and its linearization as a smooth map which is linear in $\theta$ and $\bar{\nabla}\theta$, and thus obtain a similar description for $V_2(\tilde{\varphi})-V_2(\bar{\varphi})-V_{2*}(\theta)$.

Combining these expressions for $V_1(\tilde{\varphi})$ and $V_2(\tilde{\varphi})$ we obtain
\begin{align}\label{lin-V-full}
  V(\tilde{\varphi})-V_{*}(\theta)=V(\tilde{\varphi})-V(\bar{\varphi})-V_{*}(\theta)  =F_3(\bar{\varphi},\tilde{\varphi},\theta,\bar{\nabla}\theta),
\end{align}
where $V(\bar{\varphi})=0$, $V_{*}(\theta)$ is the linearization of the vector field $V(\tilde{\varphi})$ applied to $\theta$ and $F_3$ is linear in $\theta$ and $\bar{\nabla}\theta$.  
Notice that by \eqref{lin-V}, a similar expression for $V_{2*}(\theta)$ and the linear relation between $h$ and $\theta$, we have that $V_*(\theta)$ is a linear function of $\bar{\nabla}\theta$.

Since $d\tilde{\varphi}=0$, 
we see from Cartan's formula that
\[
\mathcal{L}_{V(\tilde{\varphi})}\tilde{\varphi}
=d(V(\tilde{\varphi})\lrcorner\tilde{\varphi}).
\]
Moreover, as the exterior derivative is a linear map independent of $\theta$, and $V(\bar{\varphi})=0$, we have
\begin{align*}
D_{\bar{\varphi}}(\mathcal{L}_{V(\tilde{\varphi})}\tilde{\varphi})(\theta)&
=D_{\bar{\varphi}}\big(d(V(\tilde{\varphi})\lrcorner\tilde{\varphi})\big)(\theta)\\
&=d\big(D_{\bar{\varphi}}V(\tilde{\varphi})(\theta)\lrcorner\bar{\varphi}\big)\\
&=d(V_*(\theta)\lrcorner\bar{\varphi}).
\end{align*}
These observations together with \eqref{lin-V-full} imply that
\begin{align*}
  \mathcal{L}_{V(\tilde{\varphi})}\tilde{\varphi}-\mathcal{L}_{V(\bar{\varphi})}\bar{\varphi}-D_{\bar{\varphi}}(\mathcal{L}_{V(\tilde{\varphi})}\tilde{\varphi})(\theta)=&d(V(\tilde{\varphi})\lrcorner \tilde{\varphi}-V_{*}(\theta)\lrcorner\bar{\varphi}) \nonumber\\
   =&d\big((V(\tilde{\varphi})-V_*(\theta))\lrcorner\tilde{\varphi}+V_*(\theta)\lrcorner\theta\big)\nonumber\\
   =& d\big(F_3(\bar{\varphi},\tilde{\varphi},\theta,\bar{\nabla}\theta)\lrcorner\bar{\varphi}+V_*(\theta)\lrcorner\theta\big).
   \end{align*}
Since $V_*(\theta)$ depends linearly on $\bar{\nabla}\theta$, we deduce that
\begin{equation}   
  \mathcal{L}_{V(\tilde{\varphi})}\tilde{\varphi}-\mathcal{L}_{V(\bar{\varphi})}\bar{\varphi}-D_{\bar{\varphi}}(\mathcal{L}_{V(\tilde{\varphi})}\tilde{\varphi})(\theta)  
   =dF_4(\bar{\varphi},\tilde{\varphi},\theta,\bar{\nabla}\theta),\label{pf-linear-2}
\end{equation}
 where $F_4$ is smooth in its arguments and linear in $\theta$ and $\bar{\nabla}\theta$.
 Combining \eqref{pf-linerar-1} and \eqref{pf-linear-2} gives:
\begin{equation}\label{F-tilde}
  \tilde{F}=d\biggl(F_1(\bar{\varphi},\tilde{\varphi},\theta,\bar{\nabla}\theta)+F_4(\bar{\varphi},\tilde{\varphi},\theta,\bar{\nabla}\theta)\biggr)=dF(\bar{\varphi},\tilde{\varphi}(t),\theta(t),\bar{\nabla}\theta(t)),
\end{equation}
where $F$ is a $2$-form which is smooth in the first two arguments and linear in the last two arguments. Then \eqref{flow-theta} follows from \eqref{dt-theta}, \eqref{lin-Delta} and \eqref{F-tilde}.
\endproof

From Lemma \ref{lem-est-t-leq1} we know that $\theta(t)$ exists at least for $t\in [0,\eta]$ with $\eta\geq 1$ and
$\|\theta(t)\|_{C^k_{\bar{g}}}<\epsilon$. We next improve this estimate in a sequence of results which give estimates on $W^{m,2}$ norms
of $\theta(t)$.

In the rest of the paper, we use the same symbol $C$ to denote various positive constants which are uniformly bounded and depend at most on $\bar{\varphi}$, $\epsilon, k$ and Sobolev embedding constants.  Unless otherwise stated, all norms, inner products and integrals will be calculated with
respect to the fixed metric $\bar{g}$.

We begin with  an $L^2$ estimate, where we get exponential decay.

\begin{lem}\label{lem-estimate-1}
There exists $\bar{\epsilon}=\bar{\epsilon}(\bar{\varphi})>0$ such that if $\epsilon\in(0,\bar{\epsilon})$, then for all time $t$ for which the solution $\tilde{\varphi}(t)$ of $\eqref{deturck-flow}$ exists  and $\|\theta(t)\|_{C^k_{\bar{g}}}<\epsilon$, we have
\begin{equation}\label{L2-estimate}
  \int_M|\theta(t)|^2dv_{\bar{g}}\leq e^{-\lambda_1 t/2}\int_M|\theta(0)|^2dv_{\bar{g}},
\end{equation}
where $\lambda_1$ is the first eigenvalue of the Hodge-Laplacian
$\Delta_{\bar{\varphi}}$ on exact 3-forms, which is positive due to the Hodge decomposition theorem.
\end{lem}
\proof
We take the $\bar{g}$-inner product of \eqref{flow-theta} with $\theta(t)$ and integrate over $M$:
\begin{align}\label{proof_L2-1}
  \frac d{dt}\int_M|\theta(t)|^2dv_{\bar{g}} =& 2\int_M\theta(t)\cdot \frac \pt{\pt t}\theta(t) dv_{\bar{g}}\nonumber\\
  = & 2\int_M\theta\cdot (-\Delta_{\bar{\varphi}}\theta+dF(\bar{\varphi},\tilde{\varphi},\theta,\bar{\nabla}\theta))dv_{\bar{g}}\nonumber\\
  =&-2\int_M|d^*\theta|^2dv_{\bar{g}}+2\int_Md^*\theta\cdot Fdv_{\bar{g}}\nonumber\\
  \leq &-2\int_M|d^*\theta|^2dv_{\bar{g}}+2C\epsilon \int_M|d^*\theta||\bar{\nabla}\theta|dv_{\bar{g}}\nonumber\\
  \leq &-2\int_M|d^*\theta|^2dv_{\bar{g}}+C\epsilon \int_M(|d^*\theta|^2+|\bar{\nabla}\theta|^2)dv_{\bar{g}},
\end{align}
where $d^*=d^*_{\bar{\varphi}}$ is defined with respect to $\bar{\varphi}$ and in the first inequality of \eqref{proof_L2-1} we used the estimate $|F|\leq C|\theta||\bar{\nabla}\theta|\leq C\epsilon|\bar{\nabla}\theta|$. Recall the Weitzenb\"{o}ck formula:
\begin{equation}\label{Weitzenbo}
  \Delta_{\bar{\varphi}}\theta=-\bar{\Delta}\theta+\mathcal{R}(\theta),
\end{equation}
where $\mathcal{R}$ is some combination of curvature operators with respect to $\bar{g}$ and $\bar{\Delta}=tr\bar{\nabla}^2$ is called the rough (or connection) Laplacian.

Note that $\theta$ is closed (as it is exact). Taking the inner product with $\theta$ on both sides of \eqref{Weitzenbo} and integrating over $M$, we have
\begin{equation}\label{Weitzenbo-2}
  \int_M|d^*\theta|^2dv_{\bar{g}}=\int_M(|\bar{\nabla}\theta|^2+\theta\cdot\mathcal{R}(\theta))dv_{\bar{g}},
\end{equation}
which implies
\begin{equation}\label{proof-L2-2}
  \int_M|\bar{\nabla}\theta|^2dv_{\bar{g}}\leq \int_M(|d^*\theta|^2+C|\theta|^2)dv_{\bar{g}},
\end{equation}
where $C$ depends on the bound of the curvature tensor of $\bar{g}$. 
Since $\lambda_1$ is the first eigenvalue of $\Delta_{\bar{\varphi}}$ on exact 3-forms, we have:
\begin{equation}\label{lambda-1}
  \int_M\theta\cdot\Delta_{\bar{\varphi}}\theta dv_{\bar{g}}=\int_M|d^*\theta|^2dv_{\bar{g}}\geq \lambda_1\int_M|\theta|^2dv_{\bar{g}}.
\end{equation}
Substituting \eqref{proof-L2-2} into \eqref{proof_L2-1} and using \eqref{lambda-1}, we have
\begin{align}
  \frac d{dt}\int_M|\theta(t)|^2dv_{\bar{g}} \leq & -2(1-C\epsilon) \int_M|d^*\theta|^2dv_{\bar{g}}+C\epsilon \int_M|\theta|^2dv_{\bar{g}} \nonumber\\
  \leq & -(\lambda_1-C\epsilon)\int_M|\theta|^2dv_{\bar{g}}\nonumber\\
  \leq&-\frac{\lambda_1}2\int_M|\theta|^2dv_{\bar{g}},\label{proof-L2-3}
\end{align}
where we have chosen $\epsilon$ small such that $C\epsilon<\min\{\frac 12,\frac {\lambda_1}2\}$. Then \eqref{L2-estimate} follows from integrating \eqref{proof-L2-3} in $t$.
\endproof

We now proceed to look at higher order estimates, which will follow from the next three lemmas.

\begin{lem}\label{lem-estimate-2}
There exists  $\bar{\epsilon}=\bar{\epsilon}(\bar{\varphi})>0$ such that if $\epsilon\in(0,\bar{\epsilon})$, then for all time $t$ for which the solution $\tilde{\varphi}(t)$ of $\eqref{deturck-flow}$ exists and $\|\theta(t)\|_{C^k_{\bar{g}}}<\epsilon$, we have
\begin{equation}\label{estimate-2}
  \int_0^t\int_M|\bar{\nabla}\theta(s)|^2dv_{\bar{g}}ds\leq \left(1+\frac{C}{\lambda_1}\right)\int_M|\theta(0)|^2dv_{\bar{g}},
\end{equation}
where $\lambda_1$ is the first eigenvalue of $\Delta_{\bar{\varphi}}$
on exact 3-forms.
\end{lem}
\proof
Substituting the Weitzenb\"{o}ck formula \eqref{Weitzenbo} into the evolution equation \eqref{flow-theta}, we have
\begin{equation}\label{flow-theta-2}
  \frac{\pt}{\pt t}\theta(t)=\bar{\Delta}\theta(t)-\mathcal{R}(\theta(t))+dF.
\end{equation}
Taking the $\overline{g}$-inner product of \eqref{flow-theta-2} with $\theta(t)$ and integrating over $M$, we can use the estimates for $\mathcal{R}$ and $F$ to compute:
\begin{align}
 \frac 12 \frac d{dt}\int_M|\theta(t)|^2dv_{\bar{g}} \leq & -\int_M|\bar{\nabla}\theta|^2dv_{\bar{g}}+C\int_M|\theta|^2dv_{\bar{g}}+\int_Md^*\theta\cdot F dv_{\bar{g}}\nonumber\\
  \leq & -\int_M|\bar{\nabla}\theta|^2dv_{\bar{g}}+C\int_M|\theta|^2dv_{\bar{g}}\nonumber\\
  &\qquad+\frac 12C\epsilon\int_M(|d^*\theta|^2+|\bar{\nabla}\theta|^2)dv_{\bar{g}}\nonumber\\
  \leq&-(1-C\epsilon)\int_M|\bar{\nabla}\theta|^2dv_{\bar{g}}+C\int_M|\theta|^2dv_{\bar{g}}\nonumber\\
  \leq &-\frac 12\int_M|\bar{\nabla}\theta|^2dv_{\bar{g}}+C\int_M|\theta|^2dv_{\bar{g}},\label{pf-lem3-2-1}
\end{align}
where in the third inequality we used
\begin{equation*}
  \int_M|d^*\theta|^2dv_{\bar{g}}\leq \int_M(|\bar{\nabla}\theta|^2+C|\theta|^2)dv_{\bar{g}},
\end{equation*}
which follows from \eqref{Weitzenbo-2}. Integrating the above inequality in $t$ gives:
\begin{align*}
  \frac 12\int_0^t\int_M |\bar{\nabla}\theta(s)|^2dv_{\bar{g}}ds+\frac 12 &\int_M|\theta(t)|^2dv_{\bar{g}}\\
  &\leq  \frac 12\int_M|\theta(0)|^2dv_{\bar{g}}+C\int_0^t\int_M|\theta(s)|^2dv_{\bar{g}}ds \\
  &\leq \left(\frac 12+\frac{2C}{\lambda_1}(1-e^{-\frac{\lambda_1}2t})\right)\int_M|\theta(0)|^2dv_{\bar{g}},
\end{align*}
where in the last inequality we used the estimate
\eqref{L2-estimate}.  The result follows.
\endproof

\begin{lem}\label{lem-estimate-3}
There exists $\bar{\epsilon}=\bar{\epsilon}(\bar{\varphi})>0$ such that if $\epsilon\in(0,\bar{\epsilon})$, then for all time $t$ for which the solution $\tilde{\varphi}(t)$ of $\eqref{deturck-flow}$ exists and $\|\theta(t)\|_{C^k_{\bar{g}}}<\epsilon$, we have
\begin{align}\label{estimate-3}
   \int_0^t\int_M|\frac \pt{\pt s}\theta(s)|^2dv_{\bar{g}}ds+\int_M |\bar{\nabla} \theta&(t)|^2 dv_{\bar{g}}+\frac 12\int_0^t\int_M|\bar{\nabla}^2\theta(s)|^2dv_{\bar{g}}ds\nonumber\\
  \leq &\int_M|\bar{\nabla}\theta(0)|^2dv_{\bar{g}}+C\int_M|\theta(0)|^2dv_{\bar{g}}.
\end{align}
\end{lem}
\proof
We rewrite the equation \eqref{flow-theta-2} as
\begin{equation*}
  \frac{\pt}{\pt t}\theta(t)-\bar{\Delta}\theta(t)=-\mathcal{R}(\theta(t))+dF.
\end{equation*}
Squaring both sides with $\bar{g}$
and integrating over $M$, we can use the estimates of $\mathcal{R}$ and $F$ to see that
\begin{align}\label{eq1-est3}
  \int_M\left(\left|\frac \pt{\pt t}\theta\right|^2-2\frac \pt{\pt t}\theta\cdot \bar{\Delta}\theta+\left|\bar{\Delta}\theta\right|^2\right)dv_{\bar{g}}\leq & \int_M2(C|\theta|^2+|dF|^2)dv_{\bar{g}}.
\end{align}
By integration by parts and the Ricci identity, we compute
\begin{align}
 -\int_M2\frac \pt{\pt t}\theta\cdot \bar{\Delta}\theta dv_{\bar{g}}=& \frac d{dt}\int_M|\bar{\nabla}\theta|^2dv_{\bar{g}}, \label{eq2-est3}\\
  \int_M |\bar{\Delta}\theta|^2dv_{\bar{g}}=&\int_M|\bar{\nabla}^2\theta|^2dv_{\bar{g}}+\int_MRm(\bar{g})*\bar{\nabla}\theta*\bar{\nabla}\theta dv_{\bar{g}},\label{eq3-est3}
\end{align}
where $*$ means a contraction of tensors using  $\bar{g}$. Recall that, by the proof of Lemma \ref{lem-linearize}, we can write $F=f(\bar{\varphi},\tilde{\varphi})\theta\bar{\nabla}\theta$ for some tensor $f$ depending on $\bar{\varphi}$ and $\tilde{\varphi}$.  Moreover, $\tilde{\varphi}=\bar{\varphi}+\theta$ and $\bar{\nabla}\bar{\varphi}=0$ so $\bar{\nabla}\tilde{\varphi}=\bar{\nabla}\theta$.
Hence, differentiating the expression for $F$ yields
\begin{align*}
  |dF|\leq |\bar{\nabla} F|&\leq C\left(|\theta||\bar{\nabla}\theta|^2+|\bar{\nabla}\theta|^2+|\theta||\bar{\nabla}^2\theta|\right),
\end{align*}
where the first term on the right-hand side arises from differentiating $f$.  
Using the facts that $|\theta|<\epsilon$ and $|\bar{\nabla}\theta|<\epsilon$, we obtain the estimate
\begin{align}
  |dF|&\leq C\epsilon (|\bar{\nabla}\theta|+|\bar{\nabla}^2\theta|).\label{eq4-est3}
\end{align}
Combining \eqref{eq1-est3}-\eqref{eq4-est3} gives that
\begin{align*}
    \int_M\left|\frac \pt{\pt t}\theta\right|^2dv_{\bar{g}}+\frac d{dt}\int_M&|\bar{\nabla}\theta|^2dv_{\bar{g}}+\int_M|\bar{\nabla}^2\theta|^2dv_{\bar{g}} \\
    \leq &\int_MC\left(|\theta|^2+|\bar{\nabla}\theta|^2+ \epsilon^2|\bar{\nabla}^2\theta|^2\right)dv_{\bar{g}}.
\end{align*}
Choosing $\epsilon$ small such that $C\epsilon^2<1/2$, we have
\begin{align}
    \int_M\left|\frac \pt{\pt t}\theta\right|^2dv_{\bar{g}}+\frac d{dt}\int_M|\bar{\nabla}\theta|^2dv_{\bar{g}}+\frac 12&\int_M|\bar{\nabla}^2\theta|^2dv_{\bar{g}} \nonumber\\
   \leq &
   \int_MC\left(|\theta|^2+|\bar{\nabla}\theta|^2\right)dv_{\bar{g}}.\label{pf-lem3-3-1}
\end{align}
Integrating \eqref{pf-lem3-3-1} and using the estimates in Lemmas \ref{lem-estimate-1}-\ref{lem-estimate-2} gives the result.
\endproof

Taking covariant derivative with respect to $\bar{g}$ in \eqref{flow-theta-2}, we obtain
\begin{align}
  \frac \pt{\pt t}\bar{\nabla}\theta =& \bar{\nabla}\bar{\Delta}\theta-\mathcal{R}(\bar{\nabla}\theta)-(\bar{\nabla}\mathcal{R})\theta-\bar{\nabla} dF \nonumber\\
  =&\bar{\Delta}\bar{\nabla}\theta -\mathcal{R}_1(\bar{\nabla}\theta)-(\bar{\nabla}\mathcal{R}_1)\theta-\bar{\nabla } dF,\label{flow-nab-theta}
\end{align}
where in the second equality we used the Ricci identity
\begin{equation*}
  \bar{\nabla}\bar{\Delta}\theta-\bar{\Delta}\bar{\nabla}\theta=\bar{\nabla}Rm(\bar{g})*\theta+Rm(\bar{g})*\bar{\nabla}\theta,
\end{equation*}
and $\mathcal{R},\mathcal{R}_1$ are some combinations of curvature operators with respect to the metric $\bar{g}$.  Recall that we can write 
$F=f(\bar{\varphi},\tilde{\varphi})\theta\bar{\nabla}\theta$ and that $\bar{\nabla}\bar{\varphi}=0$ and $\bar{\nabla}\tilde{\varphi}=\bar{\nabla}\theta$.  By taking a second derivative of this expression for $F$ we deduce an estimate on $\bar{\nabla} dF$:
\begin{equation*}
  |\bar{\nabla} dF|\leq |\bar{\nabla}^2 F| \leq  C(|\theta||\bar{\nabla}\theta|^3+|\bar{\nabla}\theta|^3+
  |\bar{\nabla}^2\theta||\bar{\nabla}\theta|+|\theta||\bar{\nabla}^2\theta||\bar{\nabla}\theta|+|\theta||\bar{\nabla}^3\theta|).
  \end{equation*}
Then, using $|\theta|<\epsilon$ and $|\bar{\nabla}\theta|<\epsilon$,  we obtain
\begin{align}
  \label{nab-dF}
 |\bar{\nabla} dF| &\leq C\epsilon(|\bar{\nabla}\theta|+|\bar{\nabla}^2\theta|+|\bar{\nabla}^3\theta|).
\end{align}

\begin{lem}\label{lem-estimate-5}
There exists $\bar{\epsilon}=\bar{\epsilon}(\bar{\varphi})>0$ such that if $\epsilon\in(0,\bar{\epsilon})$, then for all time $t$ for which the solution $\tilde{\varphi}(t)$ of $\eqref{deturck-flow}$ exists and $\|\theta(t)\|_{C^k_{\bar{g}}}<\epsilon$, we have
\begin{align*}
  \int_0^t\int_M&\left|\bar{\nabla} \frac \pt{\pt s}\theta\right|^2dv_{\bar{g}}ds+\int_M|\bar{\nabla}^2\theta(t)|^2dv_{\bar{g}}+\frac 12\int_0^t\int_M|\bar{\nabla}^3\theta(s)|^2dv_{\bar{g}}ds\nonumber\\
& \leq
 \int_M|\bar{\nabla}^2\theta(0)|^2dv_{\bar{g}}+C\int_M|\bar{\nabla}\theta(0)|^2dv_{\bar{g}}+C\int_M|\theta(0)|^2dv_{\bar{g}}.
\end{align*}
\end{lem}
\proof We rewrite \eqref{flow-nab-theta} as
\begin{align}
  \frac \pt{\pt t}\bar{\nabla}\theta -\bar{\Delta}\bar{\nabla}\theta =-\mathcal{R}_1(\bar{\nabla}\theta)-(\bar{\nabla}\mathcal{R}_1)\theta-\bar{\nabla} dF.\label{flow-nab-theta-2}
\end{align}
In the usual way, squaring both sides of \eqref{flow-nab-theta-2} and integrating over $M$ gives
\begin{align}
  \int_M\left|\bar{\nabla}\frac \pt{\pt t}\theta\right|^2&dv_{\bar{g}}+\frac d{dt}\int_M|\bar{\nabla}^2\theta|^2dv_{\bar{g}}+\int_M|\bar{\Delta}\bar{\nabla}\theta|^2dv_{\bar{g}}  \nonumber\\
  &\leq C\int_M|\bar{\nabla}\theta|^2dv_{\bar{g}}+C\int_M|\theta|^2dv_{\bar{g}}+3\int_M|\bar{\nabla} dF|^2 dv_{\bar{g}}.\label{proof-est-5-1}
\end{align}
Note that curvature identities give that
\begin{equation}\label{proof-est-5-2}
\int_M|\bar{\Delta}\bar{\nabla}\theta|^2dv_{\bar{g}}=\int_M|\bar{\nabla}^3\theta|^2dv_{\bar{g}}+\int_MRm(\bar{g})*\bar{\nabla}^2\theta*\bar{\nabla}^2\theta dv_{\bar{g}}.
\end{equation}
Substituting \eqref{nab-dF} and \eqref{proof-est-5-2} into \eqref{proof-est-5-1} gives
\begin{align}
  \int_M&\left|\bar{\nabla}\frac \pt{\pt t}\theta\right|^2dv_{\bar{g}}+\frac d{dt}\int_M|\bar{\nabla}^2\theta|^2dv_{\bar{g}}+\int_M|\bar{\nabla}^3\theta|^2dv_{\bar{g}}  \nonumber\\
  &\leq C\int_M\left(|\bar{\nabla}^2\theta|^2+|\bar{\nabla}\theta|^2+|\theta|^2\right)dv_{\bar{g}}+3\int_M|\bar{\nabla} dF|^2 dv_{\bar{g}}\nonumber\\
  &\leq C\int_M\left(|\bar{\nabla}^2\theta|^2+|\bar{\nabla}\theta|^2+|\theta|^2\right)dv_{\bar{g}}+C\epsilon^2\int_M|\bar{\nabla}^3\theta|^2dv_{\bar{g}}.\label{pf-lem3-5-1}
\end{align}
Choose $\epsilon$ small such that $C\epsilon^2<1/2$. Integrating
\eqref{pf-lem3-5-1} in $t$, we get
\begin{align*}
  \int_0^t&\int_M\left|\bar{\nabla}\frac \pt{\pt s}\theta\right|^2dv_{\bar{g}}ds+\int_M|\bar{\nabla}^2\theta(t)|^2dv_{\bar{g}}+\frac 12\int_0^t\int_M|\bar{\nabla}^3\theta(s)|^2dv_{\bar{g}}ds\\
  \leq &\int_M|\bar{\nabla}^2\theta(0)|^2dv_{\bar{g}}+C\int_0^t\int_M\left(|\bar{\nabla}^2\theta(s)|^2+|\bar{\nabla}\theta(s)|^2+|\theta(s)|^2\right)dv_{\bar{g}}ds\\
  \leq &\int_M|\bar{\nabla}^2\theta(0)|^2dv_{\bar{g}}+C\int_M|\bar{\nabla}\theta(0)|^2dv_{\bar{g}}+C\int_M|\theta(0)|^2dv_{\bar{g}}
\end{align*}
as desired.  Note that we used Lemmas \ref{lem-estimate-2} and \ref{lem-estimate-3} in the derivation.
\endproof

Combining Lemmas \ref{lem-estimate-1}, \ref{lem-estimate-3} and  \ref{lem-estimate-5}, we have our Sobolev estimates:
\begin{equation}\label{pf-lem3-5-0}
  \int_{M}|\bar{\nabla}^m\theta(t)|^2dv_{\bar{g}}\leq C\|\theta(0)\|_{W^{m,2}_{\bar{g}}},\quad m=0,1,2.
\end{equation}

It is evident that we can continue repeating a similar procedure as follows. Taking $(m-1)$ covariant derivatives of \eqref{flow-theta-2} with respect to $\bar{g}$, we have
\begin{align*}
  \bar{\nabla}^{m-1}\frac{\pt}{\pt t}\theta =& \bar{\nabla}^{m-1}\bar{\Delta}\theta-\bar{\nabla}^{m-1}(\mathcal{R}(\theta))+\bar{\nabla}^{m-1}dF \\
  = &\bar{\Delta}\bar{\nabla}^{m-1}\theta+\sum_{j=0}^{m-1}\bar{\nabla}^j\mathcal{R}*\bar{\nabla}^{m-1-j}\theta+\bar{\nabla}^{m-1}dF.
\end{align*}
We rewrite the above equation as
\begin{equation*}
 \bar{\nabla}^{m-1}\frac{\pt}{\pt t}\theta -\bar{\Delta}\bar{\nabla}^{m-1}\theta=\sum_{j=0}^{m-1}\bar{\nabla}^j\mathcal{R}*\bar{\nabla}^{m-1-j}\theta+\bar{\nabla}^{m-1}dF.
\end{equation*}
Squaring both sides and integrating over $M$, we obtain
\begin{align}\label{pf-lem3-5-2}
  \int_M&\left|\bar{\nabla}^{m-1}\frac{\pt}{\pt t}\theta\right|^2dv_{\bar{g}}+\frac{d}{dt}\int_M|\bar{\nabla}^m\theta|^2dv_{\bar{g}}+\int_M|\bar{\Delta}\bar{\nabla}^{m-1}\theta|^2dv_{\bar{g}} \nonumber\\
  &\leq  C\sum_{j=0}^{m-1}\int_M|\bar{\nabla}^j\theta|^2dv_{\bar{g}}+(m+1)\int_M|\bar{\nabla}^{m-1}dF|^2dv_{\bar{g}},
\end{align}
where $C$ only depends on $\bar{\varphi}$ and $m$.  Once again using the expression 
$F=f(\bar{\varphi},\tilde{\varphi})\theta\bar{\nabla}\theta$, together with $\bar{\nabla}\bar{\varphi}=0$ and $\bar{\nabla}\tilde{\varphi}=\bar{\nabla}\theta$, we estimate $|\bar{\nabla}^{m-1}dF|\leq |\bar{\nabla}^mF|$ as follows:
\begin{align}\label{pf-lem3-5-3}
  |\bar{\nabla}^mF|   &\leq  C\sum_{j=0}^m|\bar{\nabla}^j\theta||\bar{\nabla}^{m+1-j}\theta|+C\sum_{l=1}^m|\bar{\nabla}^l\theta|\sum_{j=0}^{m-l}|\bar{\nabla}^j\theta||\bar{\nabla}^{m-l+1-j}\theta|\nonumber\\
  &\leq C\epsilon \sum_{j=[\frac m2]}^{m+1}|\bar{\nabla}^j\theta|,
\end{align}
where we need $|\bar{\nabla}^j\theta|<\epsilon$ for $0\leq j\leq [\frac m2]+1$, which will hold if $k\geq [\frac m2]+1$. Using  integration by parts and the formula for commuting covariant derivatives, we also have
\begin{align}\label{pf-lem3-5-4}
  \int_M|\bar{\Delta}\bar{\nabla}^{m-1}\theta|^2dv_{\bar{g}} = & \int_M|\nabla^{m+1}\theta|^2dv_{\bar{g}} +\int_MRm(\bar{g})*\bar{\nabla}^m\theta*\bar{\nabla}^m\theta dv_{\bar{g}}.
\end{align}
Substituting \eqref{pf-lem3-5-3} and \eqref{pf-lem3-5-4} into \eqref{pf-lem3-5-2}, we obtain
\begin{align}
  \int_M\left|\bar{\nabla}^{m-1}\frac{\pt }{\pt t}\theta\right|^2dv_{\bar{g}}&+\frac{d}{dt}\int_M|\bar{\nabla}^m\theta|^2dv_{\bar{g}}+\int_M|\bar{\nabla}^{m+1}\theta|^2dv_{\bar{g}} \nonumber\\
 & \leq  C\sum_{j=1}^{m}\int_M|\bar{\nabla}^j\theta|^2dv_{\bar{g}}+C\epsilon^2\int_M|\bar{\nabla}^{m+1}\theta|^2dv_{\bar{g}},
\end{align}
where $C=C(\bar{\varphi},m)$. Choosing $\epsilon$ such that $C\epsilon^2< 1/2$, we have
\begin{align}\label{pf-lem3-5-5}
  \int_M\left|\bar{\nabla}^{m-1}\frac{\pt}{\pt t}\theta\right|^2dv_{\bar{g}}+\frac{d}{dt}\int_M|\bar{\nabla}^m\theta|^2dv_{\bar{g}}&+\frac 12\int_M|\bar{\nabla}^{m+1}\theta|^2dv_{\bar{g}} \nonumber\\
  &\leq  C\sum_{j=1}^{m}\int_M|\bar{\nabla}^j\theta|^2dv_{\bar{g}}.
\end{align}
Starting from \eqref{pf-lem3-5-0}, an induction argument applied to \eqref{pf-lem3-5-5} gives that
\begin{align}\label{pf-lem3-5-6}
  \int_0^t\int_M\left|\bar{\nabla}^{m-1}\frac{\pt}{\pt s}\theta\right|^2&dv_{\bar{g}}ds+\int_M|\bar{\nabla}^m\theta(t)|^2dv_{\bar{g}}+\frac 12\int_0^t\int_M|\bar{\nabla}^{m+1}\theta(s)|^2dv_{\bar{g}}ds \nonumber\\
  &\leq  \int_M|\bar{\nabla}^m\theta(0)|^2dv_{\bar{g}}+C\sum_{j=1}^{m-1}\int_M|\bar{\nabla}^j\theta(0)|^2dv_{\bar{g}}.
\end{align}
So we have
\begin{equation}\label{W^2m-estimate}
  \int_{M}|\bar{\nabla}^m\theta(t)|^2dv_{\bar{g}}\leq C\|\theta(0)\|_{W^{m,2}_{\bar{g}}}
\end{equation}
whenever $[\frac{m}{2}]\leq k-1$.  

Recall that as $\text{dim}\,M=7$ if $m-k> \frac 72$,
then  the Sobolev embedding theorem says that $W^{m,2}_{\bar{g}}$ embeds continuously in $C^{k}_{\bar{g}}$ by inclusion.
  Since $m=k+4$ satisfies $k+\frac{7}{2}< m$ trivially and $[\frac{m}{2}]\leq k-1$ for all $k\geq 5$ (since $k\geq 5$ implies $[\frac{k}{2}]\leq k-3$), we can use
 \eqref{W^2m-estimate} to obtain the main estimate of this section.

\begin{lem}\label{lem-C^k}
For any $k\geq 5$ 
there exists $\bar{\epsilon}>0$ such that if $\epsilon\in(0,\bar{\epsilon})$, then for all time $t$ for which the solution $\tilde{\varphi}(t)$ of $\eqref{deturck-flow}$ exists  and $\|\theta(t)\|_{C^k_{\bar{g}}}<\epsilon$, we have
\begin{equation}\label{C^k-est-1}
  \|\theta(t)\|_{C^k_{\bar{g}}}\leq C_1\|\theta(0)\|_{W^{k+4,2}_{\bar{g}}}
  \leq C_2\|\theta(0)\|_{C^{k+4}_{\bar{g}}},
\end{equation}
where $C_1,C_2$ are constants depending only on $\bar{\varphi},\epsilon,k$ and Sobolev embedding constants.
\end{lem}

\section{Proof of Theorem \ref{thm-main}}\label{sec-proof of thm}

In this section, we prove the main theorem of this paper, Theorem \ref{thm-main}. First, we show the dynamical stability of torsion-free $\GG_2$ structures under the Laplacian--DeTurck flow \eqref{deturck-flow} for closed $\GG_2$ structures.
\begin{thm}\label{thm-deturkflow}
Let $(M,\bar{\varphi})$ be a compact $\GG_2$ manifold, let $k\geq 5$ and let $\bar{\epsilon}=\bar{\epsilon}(\bar{\varphi})>0$ be such that Lemmas \ref{lem-linearize}-\ref{lem-C^k} hold.  Let $\epsilon\in(0,\bar{\epsilon})$ and let $C_1$ be the constant given in \eqref{C^k-est-1}.

There exists $\delta=\delta(M,\bar{\varphi},k,\epsilon)>0$ such that if $\varphi_0\in [\bar{\varphi}]_+$ satisfies $$\|\varphi_0-\bar{\varphi}\|_{C^{k+4}_{\bar{g}}}<\delta\quad\text{and}\quad \|\varphi_0-\bar{\varphi}\|_{W^{k+4,2}_{\bar{g}}}\leq \epsilon/{2C_1},$$
then the solution $\tilde{\varphi}(t)$ of the Laplacian--DeTurck flow \eqref{deturck-flow} starting at $\varphi_0$ exists for all time $t\in [0,\infty)$, satisfies $\|\tilde{\varphi}(t)-\bar{\varphi}\|_{C^{k}_{\bar{g}}}<\epsilon$ for all $t\in [0,\infty)$ and $\tilde{\varphi}(t)\ra \bar{\varphi}$ exponentially in $C^k_{\bar{g}}$ as $t\ra \infty$.
\end{thm}
\proof
By Lemma \ref{lem-est-t-leq1} there exists $\delta=\delta(M,\bar{\varphi},k,\epsilon)>0$
such if  the initial value $\varphi_0\in [\bar{\varphi}]_+$ satisfies $\|\varphi_0-\bar{\varphi}\|_{C^{k+4}_{\bar{g}}}<\delta$, then the solution $\tilde{\varphi}(t)$ of the flow \eqref{deturck-flow} exists at least
for time $t\in[0,1]$.  Let $\theta(t)=\tilde{\varphi}(t)-\bar{\varphi}$ as in the previous section.

Let $\eta$ be the maximum time such that $\tilde{\varphi}(t)$ exists and $\|\theta(t)\|_{C^{k}_{\bar{g}}}<\epsilon$ for all $t\in [0,\eta)$. By Lemma \ref{lem-est-t-leq1}, we know that $\eta\geq 1$. 
 Suppose towards a contradiction that $\eta<\infty$. From the estimates \eqref{C^k-est-1}, we have that
\begin{equation}
   \|\theta(t)\|_{C^k_{\bar{g}}}\leq C_1\|\theta(0)\|_{W^{k+4,2}_{\bar{g}}}\leq C_1\frac{\epsilon}{2C_1}= \epsilon/2,
   \end{equation}
for all $t\in [0,\eta)$. This contradicts the maximality in the definition of $\eta$. Thus $\eta=\infty$ and
\begin{equation}
   \|\theta(t)\|_{C^k_{\bar{g}}}<\epsilon \quad \forall t\in [0,\infty).\end{equation}

To complete the proof, we need to demonstrate the exponential  convergence of $\tilde{\varphi}(t)\ra \bar{\varphi}$  as $t\ra \infty$ in $C^k_{\bar{g}}$. First, by Lemma \ref{lem-estimate-1}, $\varphi(t)\ra \bar{\varphi}$ exponentially in $L^2_{\bar{g}}$, so $\theta(t)\to 0$
in $L^2_{\bar{g}}$. 
Now, from \eqref{pf-lem3-2-1} in the proof of Lemma \ref{lem-estimate-2},
\begin{align*}
 \frac 12 \frac d{dt}\int_M|\theta(t)|^2dv_{\bar{g}}
 &\leq-\frac 12\int_M|\bar{\nabla}\theta|^2dv_{\bar{g}}
 +C\int_M|\theta|^2dv_{\bar{g}}.
\end{align*}
Integrating this equation in time $t$, using estimate \eqref{L2-estimate} and the fact that $\theta(t)\to 0$ in $L^2$ we get
\begin{align}
 \frac 12\int_t^{\infty}\int_M|\bar{\nabla}\theta(s)|^2dv_{\bar{g}}ds
 &\leq  \frac 12 \int_M|\theta(t)|^2dv_{\bar{g}} 
 +C\int_t^{\infty}\int_M|\theta(s)|^2dv_{\bar{g}}ds \nonumber \\
  &\leq   \left(\frac 12+\frac{2C}{\lambda_1}\right)e^{-\frac{\lambda_1 t}2} \int_M|\theta(0)|^2dv_{\bar{g}}.\label{pf-converg-0}
\end{align}
Recall \eqref{pf-lem3-3-1} in the proof of Lemma \ref{lem-estimate-3}:\begin{align}
    \int_M\left|\frac \pt{\pt t}\theta\right|^2dv_{\bar{g}}+\frac d{dt}\int_M|\bar{\nabla}\theta|^2dv_{\bar{g}}+&\frac 12\int_M|\bar{\nabla}^2\theta|^2dv_{\bar{g}} \nonumber \\
   &\leq \int_MC\left(|\theta|^2+|\bar{\nabla}\theta|^2\right)dv_{\bar{g}}.\label{pf-converg-1}
\end{align}
Let $\xi(s)$ be a cut-off function with $0\leq\xi(s)\leq 1$  satisfying
$\xi(s)=0$ for $s\leq t-\frac 12$, $|\xi'(s)|\leq C$ for $s\in [t-\frac 12,t]$ and $\xi(s)=1$ for $s\geq t$. Then from \eqref{pf-converg-1}, we have
\begin{align*}
  \frac d{ds}\bigg(\xi(s)\int_M|\bar{\nabla}&\theta(s)|^2dv_{\bar{g}}\bigg)\\
  &=  \xi'(s)\int_M|\bar{\nabla}\theta(s)|^2dv_{\bar{g}}+\xi(s)\frac d{ds}\int_M|\bar{\nabla}\theta(s)|^2dv_{\bar{g}} \displaybreak[0]\\
  &\leq  \left(\xi'(s)+C\xi(s)\right)\int_M|\bar{\nabla}\theta(s)|^2dv_{\bar{g}}+C\xi(s)\int_M|\theta(s)|^2dv_{\bar{g}}\\
  &\leq  C\int_M\left(|\bar{\nabla}\theta(s)|^2+|\theta(s)|^2\right)dv_{\bar{g}}.
\end{align*}
Integrating the above inequality from $t-\frac 12$ to $t$ and using \eqref{L2-estimate} and \eqref{pf-converg-0} we get
\begin{align}
  \int_{M}|\bar{\nabla}\theta(t)|^2dv_{\bar{g}} &\leq C\int_{t-\frac 12}^t\int_M\left(|\bar{\nabla}\theta(s)|^2+|\theta(s)|^2\right)dv_{\bar{g}} \nonumber\\
  &\leq  C\int_{t-\frac 12}^{\infty}\int_M\left(|\bar{\nabla}\theta(s)|^2+|\theta(s)|^2\right)dv_{\bar{g}} \nonumber\\
  &\leq Ce^{-\frac{\lambda_1 t}2} \int_M|\theta(0)|^2dv_{\bar{g}},
\end{align}
where $C$ is a uniform constant. 
Hence $\bar{\nabla}\theta\to 0$ in $L^2_{\bar{g}}$ as $t\ra\infty$.  Using this fact, Lemma \ref{lem-estimate-1} and \eqref{pf-converg-1}  imply that
\begin{align*}
  \int_t^{\infty}\int_M|\bar{\nabla}^2&\theta(s)|^2dv_{\bar{g}}ds\\
  &\leq  2 \int_{M}|\bar{\nabla}\theta(t)|^2dv_{\bar{g}} 
   +2C\int_t^{\infty}\int_M\left(|\bar{\nabla}\theta(s)|^2+|\theta(s)|^2\right)dv_{\bar{g}}ds\\
   &\leq Ce^{-\frac{\lambda_1 t}2} \int_M|\theta(0)|^2dv_{\bar{g}}.
\end{align*}
Using \eqref{pf-lem3-5-1} and a cut-off argument in a similar way as above, we have
\begin{align*}
 \int_{M}|\bar{\nabla}^2\theta(t)|^2dv_{\bar{g}} \leq& C\int_{t-\frac 12}^t\int_M\left(|\bar{\nabla}^2\theta(s)|^2+|\bar{\nabla}\theta(s)|^2+|\theta(s)|^2\right)dv_{\bar{g}} \nonumber\\
  \leq & C\int_{t-\frac 12}^{\infty}\int_M\left(|\bar{\nabla}^2\theta(s)|^2+|\bar{\nabla}\theta(s)|^2+|\theta(s)|^2\right)dv_{\bar{g}} \nonumber\\
  \leq &Ce^{-\frac{\lambda_1 t}2} \int_M|\theta(0)|^2dv_{\bar{g}},
\end{align*}
By repeating this procedure and using the induction part in \S \ref{sec-estimate}, we can obtain
\begin{align}
  \int_{M}|\bar{\nabla}^m\theta(t)|^2dv_{\bar{g}} \leq & Ce^{-\frac{\lambda_1 t}2} \int_M|\theta(0)|^2dv_{\bar{g}}
\end{align}
for all $m$ so that $[\frac{m}{2}]\leq k-1$, in particular for $m\leq k+4$ since $k\geq 5$. Then the Sobolev embedding theorem gives that
\begin{equation}\label{pf-converg-2}
  \|\tilde{\varphi}(t)-\bar{\varphi}\|_{C^k_{\bar{g}}}=\|\theta(t)\|_{C^k_{\bar{g}}}\leq Ce^{-\frac{\lambda_1 t}2}
\end{equation}
as required.
\endproof

 Theorem \ref{thm-deturkflow} shows that the solution $\tilde{\varphi}(t)$ of the Laplacian--DeTurck flow \eqref{deturck-flow} converges to the torsion free $\GG_2$ structure $\bar{\varphi}$ exponentially in the $C^k_{\bar{g}}$-norm, if the initial value $\varphi_0$ is sufficiently close to $\bar{\varphi}$. Let $\varphi(t)$ be the solution to the Laplacian flow \eqref{Lap-flow-def} with the same initial value $\varphi_0$. Then $\varphi(t)=\phi_t^*\tilde{\varphi}(t)$ and the associated metric $g(t)=\phi_t^*\tilde{g}(t)$, where $\phi_t$ is a family of diffeomorphisms satisfying
\begin{equation}\label{phi-t}
  \left\{\begin{array}{ccl}
           \dfrac \pt{\pt t}\phi_t(x) & = & -V(\tilde{\varphi}(t))|_{\phi_t(x)}\\
           \phi_0 & = & id
         \end{array}\right.
\end{equation}
for $x\in M$, where the vector field $V(\tilde{\varphi}(t))$ is defined in \eqref{vector-V}. Note that $V(\tilde{\varphi}(t))$ involves the first order covariant derivative of $\tilde{\varphi}(t)$. Since $\tilde{\varphi}(t)$ exists for all time $[0,\infty)$ and converges to $\bar{\varphi}$ exponentially in $C^k_{\bar{g}}$-norm where $k\geq 3$, the vector field $V(\tilde{\varphi}(t))$ also exists for all time $[0,\infty)$ and converges to $V(\bar{\varphi})=0$ in the $C^{k-1}_{\bar{g}}$-norm exponentially. Then the solution $\phi_t$ of \eqref{phi-t} exists for all $t\in[0,\infty)$ (cf.~\cite[Lemma 3.15]{Chow-Knopf}). Moreover, $\phi_t$ converges to a limit map $\phi_{\infty}$ in $C^{k-1}_{\bar{g}}$-norm as $t\ra\infty$.

Recalling that $\tilde{\varphi}(t)-\bar{\varphi}=\theta(t)$ and $\bar{\varphi}$ is torsion-free, we have by \eqref{pf-linerar-1} that
\[
\Delta_{\tilde{\varphi}(t)}\tilde{\varphi}(t)=\Delta_{\tilde{\varphi}(t)}\tilde{\varphi}(t)-\Delta_{\bar{\varphi}}\bar{\varphi}=L(\theta(t))+d F(\theta(t),\bar{\nabla}\theta(t)),
\]
where $L$ is a linear second order differential operator and $F$ is dominated by $\theta\bar{\nabla}\theta$ when $|\theta|$ is sufficiently small.  Therefore, since the exterior derivative is a linear first order differential operator, we have 
\begin{gather*}
\|L(\theta(t))\|_{C^{k-2}_{\bar{g}}}\leq C\|\theta(t)\|_{C^k_{\bar{g}}},\\
\|d F(\theta(t),\bar{\nabla}\theta(t))\|_{C^{k-2}_{\bar{g}}}\leq 
C\| F(\theta(t),\bar{\nabla}\theta(t))\|_{C^{k-1}_{\bar{g}}}
\leq C\|\theta(t)\|_{C^k_{\bar{g}}}^2
\leq C\|\theta(t)\|_{C^{k}_{\bar{g}}},
\end{gather*}
when $\|\theta(t)\|_{C^k_{\bar{g}}}$ is sufficiently small.
 Combining these estimates 
with \eqref{pf-converg-2} we obtain
\begin{align*}
  \|\Delta_{\tilde{\varphi}(t)}\tilde{\varphi}(t)\|_{C^{k-2}_{\bar{g}}}&
  \leq  C\|\tilde{\varphi}(t)-\bar{\varphi}\|_{C^k_{\bar{g}}}\leq Ce^{-\frac{\lambda_1 t}2}.
\end{align*}
Since $\|\tilde{\varphi}(t)-\bar{\varphi}\|_{C^{k}_{\bar{g}}}<\epsilon$ for all time $t\in [0,\infty)$, by making $\epsilon$ smaller if necessary we can assume that the $C^k$-norms of $\bar{g}$ and $\tilde{g}(t)$ differ at most by a factor $2$ for all $t\in[0,\infty)$. Then we have
\begin{align}\label{est-Delta-sigma}
  \|\Delta_{\tilde{\varphi}(t)}\tilde{\varphi}(t)\|_{C^{k-2}_{\tilde{g}(t)}}\leq & 2 \|\Delta_{\tilde{\varphi}(t)}\tilde{\varphi}(t)\|_{C^{k-2}_{\bar{g}}}\leq 2Ce^{-\frac{\lambda_1 t}2}
\end{align}
By diffeomorphism invariance, \eqref{est-Delta-sigma} implies that
\begin{align}\label{C^0-hat-sigma}
  \|\Delta_{\varphi(t)}\varphi(t)\|_{C^{k-2}_{g(t)}}\leq &  2Ce^{-\frac{\lambda_1 t}2}.
\end{align}

We know that for closed $\GG_2$ structures,
  $\Delta_{\varphi}\varphi=i_{\varphi}(h)$,
where $h$ is a symmetric $2$-tensor satisfying \eqref{hodge-Lap-varp-3}. By  \cite[Proposition 2.9]{Kar},
\begin{equation}\label{h-est}
  |\Delta_{\varphi}\varphi|^2_g=|i_{\varphi}(h)|^2_g=(tr_g(h))^2+2h_i^kh_k^i.
\end{equation}
From \eqref{C^0-hat-sigma}-\eqref{h-est}, we have
\begin{equation}\label{h-est2}
  \|h(t)\|_{C^{k-2}_{g(t)}}\leq \frac{C}{2}e^{-\frac{\lambda_1 t}2}.
\end{equation}
Recall that under the Laplacian flow the associated metric $g(t)$ of $\varphi(t)$ evolves by
\begin{equation*}
  \frac{\pt}{\pt t}g(t)=2h(t).
\end{equation*}
From \eqref{h-est2},
\begin{equation}\label{h-est3}
  \int_0^{\infty}\|h(t)\|_{C^{k-2}_{g(t)}}dt\leq \frac{C}{2}\int_0^{\infty}e^{-\frac{\lambda_1 t}2}dt=\frac{C}{\lambda_1}.
\end{equation}
We now observe that, for any tangent vector $X\neq 0$ on $M$, we have 
  \begin{align*}
    \left|\frac{\partial}{\partial t}g_t(X,X)\right| = 2|h_t(X,X)|\leq 2\|h_t\|_{C^0_{g(t)}} g_t(X,X), 
  \end{align*}
  which implies that
   \begin{align*}
    \left|\frac{\partial}{\partial t}\ln g_t(X,X)\right| \leq 2\|h_t\|_{C^0_{g(t)}} . 
  \end{align*}
  Then, by \eqref{h-est3}, we have 
  \begin{equation*}
    \left|\ln \frac{g_t(X,X)}{g_0(X,X)}\right|\leq \int_0^t2\|h_t\|_{C^0_{g(t)}}\leq \int_0^{\infty}2\|h_t\|_{C^0_{g(t)}}\leq \frac{2C}{\lambda_1}.
  \end{equation*}
Hence, for any $t\in [0,\infty)$,
\begin{equation}\label{eqiv-gt}
  e^{-\frac{2C}{\lambda_1}}g(0)\leq g(t)\leq e^{\frac{2C}{\lambda_1}}g(0).
\end{equation}
Note that $g(0)$ is the associated metric of $\varphi_0$ (which is very close to $\bar{\varphi}$), so $g(0)$ is uniformly equivalent to $\bar{g}$, the associated metric of $\bar{\varphi}$. Therefore, from \eqref{eqiv-gt} we have $g(t)$ is uniformly equivalent to $\bar{g}$ for all $t\in [0,\infty)$.

Then we deduce from \eqref{C^0-hat-sigma} that
\begin{align}\label{decay-lapl}
  \|\Delta_{\varphi(t)}\varphi(t)\|_{C^{0}_{\bar{g}}}\leq &  C\|\Delta_{\varphi(t)}\varphi(t)\|_{C^{0}_{g(t)}}\leq   Ce^{-\frac{\lambda_1 t}2},
\end{align}
and therefore $\|\frac {d}{dt}\varphi(t)\|_{C^0_{\bar{g}}}\leq Ce^{-\frac{\lambda_1 t}2}$. This tells us there exists a limit $\GG_2$ structure $\varphi_{\infty}$ such that
\begin{align}\label{c^0-converg}
  \|\varphi(t)-\varphi_{\infty}\|_{C^0_{\bar{g}}}\leq &\int_t^{\infty}\|\frac {d}{ds}\varphi(s)\|_{C^0_{\bar{g}}}ds\nonumber\\
  = &\int_t^{\infty} \|\Delta_{\varphi(s)}\varphi(s)\|_{C^{0}_{\bar{g}}}ds\nonumber\\
  \leq &\frac{2C}{\lambda_1}e^{-\frac{\lambda_1t}2}.
\end{align}
Note that $\varphi_{\infty}$ is a positive $3$-form because we can see as in \cite{Lotay-Wei} that $\varphi_{\infty}$ defines a metric $g_{\infty}$ which is the limit of the metric $g(t)$. The limit $\varphi_{\infty}$ is torsion-free as $M$ is compact and we have $\Delta_{\varphi_{\infty}}\varphi_{\infty}=0$ from \eqref{decay-lapl}.

We claim that the limit $\varphi_{\infty}=\phi_{\infty}^*\bar{\varphi}$. In fact, for any time-independent vector field $v$ on $M$, we see from \eqref{phi-t} that the differential $\phi_{t*}$ of $\phi_t$ satisfies
\begin{align*}
  \frac 12\frac d{dt}|\phi_{t*}(v)|_{\bar{g}}^2 =& \bar{g}(\phi_{t*}(v),\frac d{dt}\phi_{t*}(v)) \\
  = & \bar{g}(\phi_{t*}(v),\nabla_{\phi_{t*}(\pt_t)}\phi_{t*}(v))\\
   = & \bar{g}(\phi_{t*}(v),\nabla_{\phi_{t*}(v)}\phi_{t*}(\pt_t))\\
   = & \bar{g}(\phi_{t*}(v),-\nabla_{\phi_{t*}(v)}V(\tilde{\varphi}(t))).
\end{align*}
We deduce that 
\begin{equation}\label{dphi-t2}
  -\|\nabla V(\tilde{\varphi}(t)))\|_{C^0_{\bar{g}}}\leq \frac d{dt}\ln |\phi_{t*}(v)|_{\bar{g}}\leq \|\nabla V(\tilde{\varphi}(t)))\|_{C^0_{\bar{g}}}
\end{equation}
Since $V(\tilde{\varphi}(t))$ converges to zero exponentially in $C^{k-1}$-norm, then
\begin{equation*}
  \int_0^t\|\nabla V(\tilde{\varphi}(s)))\|_{C^0_{\bar{g}}}ds\leq \int_0^tCe^{-\frac{\lambda_1 s}2}ds=\frac{2C}{\lambda_1}(1-e^{-\frac{\lambda_1t}2})
\end{equation*}
Note that $\phi_0$ is the identity map, which means that $|\phi_{0*}(v)|_{\bar{g}}=|v|_{\bar{g}}$. Integrating the inequalities in \eqref{dphi-t2}, we have
\begin{align*}
 |\phi_{t*}(v)|_{\bar{g}} \leq & |\phi_{0*}(v)|_{\bar{g}}e^{\int_0^{t}\|\nabla V(\tilde{\varphi}(s)))\|_{C^0_{\bar{g}}}ds}\leq |v|_{\bar{g}}e^{\frac{2C}{\lambda_1}}
\end{align*}
and
\begin{align*}
 |\phi_{t*}(v)|_{\bar{g}} \geq & |\phi_{0*}(v)|_{\bar{g}}e^{-\int_0^{t}\|\nabla V(\tilde{\varphi}(s)))\|_{C^0_{\bar{g}}}ds}  \geq |v|_{\bar{g}}e^{-\frac{2C}{\lambda_1}}
\end{align*}
So the differential $\phi_{t*}$ of the map $\phi_t$ is non-degenerate everywhere on $M$ and the norm of $\phi_{t*}$ is uniformly bounded independently of $t$. This estimate thus also descends to the limit map $\phi_{\infty}$. The inverse function theorem then implies that $\phi_{\infty}$ is a local diffeomorphism.  Since $\phi_0=id$ the identity map, and each $\phi_t$ is a diffeomorphism which is isotopic to the identity map, we deduce that $\phi_\infty$ is a surjective local diffeomorphism homotopic to the identity.  As $M$ is compact, $\phi_{\infty}$ is a covering map.  Moreover, $\phi_{\infty}$ is homotopic to the identity so it must have degree one and thus is injective.  We deduce that the limit map $\phi_{\infty}$ is also a diffeomorphism. Then
\begin{align*}
  |\varphi_{\infty}-\phi_{\infty}^*\bar{\varphi}| \leq & \lim_{t\ra\infty}\left( |\varphi_{\infty}-\phi_t^*\tilde{\varphi}(t)|+|\phi_t^*(\tilde{\varphi}(t)-\bar{\varphi})|+|(\phi_t^*-\phi_{\infty}^*)\bar{\varphi}|\right)\\
  \leq & \lim_{t\ra\infty}\left( |\varphi_{\infty}-\varphi(t)|+C|(\tilde{\varphi}(t)-\bar{\varphi})|+C|\phi_t^*-\phi_{\infty}^*|\right)=0.
\end{align*}
We conclude that the Laplacian flow $\varphi(t)$ converges to $\varphi_{\infty}=\phi_{\infty}^*\bar{\varphi}$, which lies in the diffeomorphism orbit of $\bar{\varphi}$. Moreover, the convergence is exponentially in the $C^0$-norm.

We next show that  $\varphi(t)$ converges to $\varphi_{\infty}=\phi_{\infty}^*\bar{\varphi}$ in any $C^k$-norm. Since $\tilde{\varphi}(t)$ converges to $\bar{\varphi}$ exponentially in $C^k$-norm $(k\geq 5)$, we have
\begin{equation*}
  |Rm(\tilde{g}(t))-Rm(\bar{g})|_{\bar{g}}\leq Ce^{-\frac{\lambda_1t}2},\quad |\tilde{\nabla} \tilde{T}(x,t)|_{\bar{g}}\leq Ce^{-\frac{\lambda_1t}2}.
\end{equation*}
We also note that for $t$ sufficiently large, $\tilde{g}(t)$ is uniformly equivalent to $\bar{g}$, therefore we have
\begin{equation*}
  |Rm(\tilde{g}(t))|_{\tilde{g}(t)}+|\tilde{\nabla} \tilde{T}(x,t)|_{\tilde{g}(t)}\leq C,
\end{equation*}
for large $t$, i.e.~for $t\geq t_0$, which gives us
\begin{equation}
  |Rm({g}(t))|_{{g}(t)}+|{\nabla} {T}(x,t)|_{{g}(t)}\leq C,
\end{equation}
for $t\geq t_0$ by diffeomorphism invariance. By Shi-type derivative estimates in \cite{Lotay-Wei}, we get uniform bounds on the covariant derivatives of $Rm$ and $\nabla T$,
\begin{equation}\label{Rm-T-bd}
  |\nabla^kRm({g}(t))|_{{g}(t)}+|{\nabla}^{k+1} {T}(x,t)|_{{g}(t)}\leq C(k),
\end{equation}
for $t\geq t_0+1$ and all $k\geq 1$. The injectivity radius of $g(t)$ also satisfies a uniform lower bound $inj(g(t))\geq \delta>0$ since $g(t)$ is uniformly equivalent to $\bar{g}$ (see the proof of \cite[Theorem 8.1]{Lotay-Wei}). 
Then from the compactness theorem for $\GG_2$ structures in \cite{Lotay-Wei}, there is a sequence of times $t_i\ra \infty$ such that $(M,\varphi(t_i))$ converges smoothly to $(M_{\infty},\hat{\varphi}_{\infty})$. Note that, since the metrics $g(t)$ are uniformly equivalent to $\bar{g}$, the diameters of the manifolds $(M,g(t))$ are uniformly bounded, and thus the compactness theory implies that the limit  $M_{\infty}$ is diffeomorphic to $M$.  We can then pull back $\hat{\varphi}_{\infty}$ to $M$ and, by abuse of notation, view it as being defined on $M$. Therefore the closed $\GG_2$ structures $\varphi(t_i)$ converge smoothly to a closed $\GG_2$ structure $\hat{\varphi}_{\infty}$ on $M$ as $t_i\ra\infty$. Since we already know that $\varphi(t)$ converges to $\varphi_{\infty}$ continuously (i.e.~in the $C^0_{\bar{g}}$-norm) as $t\ra\infty$, the  uniqueness of limits implies that the limit $\GG_2$ structure $\hat{\varphi}_{\infty}=\varphi_{\infty}$. 

Moreover, we have this smooth convergence without passing to a subsequence, in the sense that $\varphi(t)\ra\varphi_{\infty}$ smoothly as $t\ra\infty$. 
For a contradication we suppose not, so that there exists an integer $k$, $\epsilon>0$ and a sequence of times $t_l\ra\infty$ such that, for all $l$,
 \begin{equation*}
   \|\varphi(t_l)-\varphi_{\infty}\|_{C^k_{\bar{g}}}\geq \epsilon.
 \end{equation*}
Since $\varphi(t_l)$ satisfies  \eqref{Rm-T-bd} and the injectivity radius estimate, we can apply the compactness theorem for $\GG_2$ structures in \cite{Lotay-Wei} again to obtain a subsequence $t_{l(j)}$ which converges in $C^k_{\bar{g}}$ to a limit $\varphi'_{\infty}$ on $M$ with
 \begin{equation*}
   \|\varphi'_{\infty}-\varphi_{\infty}\|_{C^k_{\bar{g}}}\geq \epsilon.
 \end{equation*}
Thus $\varphi'_{\infty}\neq \varphi_{\infty}$, which contradicts the fact that $\varphi(t_{l(j)})$ converges to $\varphi_{\infty}$ continuously.  

We conclude that $\varphi(t)$ converges to $\varphi_{\infty}=\phi_{\infty}^*\bar{\varphi}$ smoothly as $t\ra \infty$. Since $\phi_0$ is the identity map, we have $\varphi_{\infty}\in \Diff^{0}(M)\cdot\bar{\varphi}$. This completes the proof of Theorem \ref{thm-main}.

\bibliographystyle{Plain}

\end{document}